\documentclass[final]{siamart1116}

\usepackage{lipsum}
\usepackage{amsfonts}
\usepackage{graphicx}
\usepackage{epstopdf}
\usepackage{algorithmic}

\ifpdf
  \DeclareGraphicsExtensions{.eps,.pdf,.png,.jpg}
\else
  \DeclareGraphicsExtensions{.eps}
\fi

\usepackage{changepage}
\usepackage[caption=false]{subfig}

\newcommand{\jump}[1]{\mbox{$[\![ #1 ]\!]$}}
\newcommand{\avg}[1]{\mbox{$\{\!\>\!\>\!\!\!\{ #1 \}\!\>\!\>\!\!\!\}$}}

\newcommand{\dx}{\,\text{d}x}
\newcommand{\dd}[1]{\frac{\text{d}}{\text{d} #1}}
\newcommand{\ddd}[1]{\frac{\text{d}^2}{\text{d} #1 ^2}}
\newcommand{\limeps}{\text{{\raisebox{0.2ex}{\scalebox{0.8}{$\lim\limits_{\epsilon\rightarrow 0}$}}}}}
\renewcommand{\d}[1]{\,\text{d}#1}

\makeatletter
\newcommand*\rel@kern[1]{\kern#1\dimexpr\macc@kerna}
\newcommand*\widebar[1]{%
  \begingroup
  \def\mathaccent##1##2{%
    \rel@kern{1.2}%
    \overline{\rel@kern{-1.2}\macc@nucleus\rel@kern{-0.2}}%
    \rel@kern{-0.2}%
  }%
  \macc@depth\@ne
  \let\math@bgroup\@empty \let\math@egroup\macc@set@skewchar
  \mathsurround\z@ \frozen@everymath{\mathgroup\macc@group\relax}%
  \macc@set@skewchar\relax
  \let\mathaccentV\macc@nested@a
  \macc@nested@a\relax111{#1}%
  \endgroup
}
\makeatother

\numberwithin{theorem}{section}

\newcommand{\TheTitle}{Residual-based variational multiscale modeling in a discontinuous Galerkin framework} 
\newcommand{\TheTitleHeader}{Residual-based VMS modeling in a DG framework} 
\newcommand{\TheAuthors}{S.K.F. Stoter, S.R. Turteltaub, S.J. Hulshoff, D. Schillinger}

\headers{\TheTitleHeader}{\TheAuthors}

\title{{\TheTitle}\thanks{Submitted to the editors on 09-18-2017.
\funding{This work was funded by the National Science Foundation through the NSF CAREER Award No. 1651577.}}}

\author{
  Stein K.F. Stoter\thanks{Department of Civil, Environmental, and Geo- Engineering, University of Minnesota, USA
    (\email{Stote031@umn.edu}, \email{Dominik@umn.edu}).}
  \and
  Sergio R. Turteltaub\thanks{Faculty of Aerospace Engineering, Delft University of Technology, The Netherlands
  	(\email{S.R.Turteltaub@tudelft.nl}, \email{S.J.Hulshoff@tudelft.nl}).}
  \and
  Steven J. Hulshoff\footnotemark[3]
  \and
  Dominik Schillinger\footnotemark[2]
}

\usepackage{amsopn}

\ifpdf
\hypersetup{
  pdftitle={\TheTitle},
  pdfauthor={\TheAuthors}
}
\fi

\usepackage{multirow}
\newcommand{\specialcell}[2][l]{%
  \begin{tabular}[l]{@{}#1@{}}#2\end{tabular}}

\begin{document}

\maketitle

\begin{abstract}
We develop the general form of the variational multiscale method in a discontinuous Galerkin framework. Our method is based on the decomposition of the true solution into discontinuous coarse-scale and discontinuous fine-scale parts.
The obtained coarse-scale weak formulation includes two types of fine-scale contributions. The first type corresponds to 
a fine-scale volumetric term, which we formulate in terms of a residual-based model that also takes into account fine-scale effects at element interfaces. The second type consists of independent fine-scale terms at element interfaces, which we formulate in terms of a new fine-scale ``interface model''.
We demonstrate for the one-dimensional Poisson problem that existing discontinuous Galerkin formulations, such as the interior penalty method, can be rederived by choosing particular fine-scale interface models. The multiscale formulation thus opens the door for a new perspective on discontinuous Galerkin methods and their numerical properties. 
This is demonstrated for the one-dimensional advection-diffusion problem, where we show that upwind numerical fluxes can be interpreted as an ad hoc remedy for missing volumetric fine-scale terms.
\end{abstract}




\section{Introduction}

The variational multiscale (VMS) method was proposed by Hughes and collaborators in the 1990s \cite{Hughes1995,Hughes1998,Hughes1996} as a strategy for capturing the subgrid-scale behavior of discrete solutions of partial differential equations (PDEs) in the variational form. 
Already in the first work on the VMS method, it was recognized that stabilized formulations could be obtained by consistent incorporation of fine-scale effects \cite{Hughes1995}. Many classical stabilization techniques for advection-type problems, such as SUPG \cite{Brooks1982}, GLS \cite{HUGHES1989173} and PSPG  \cite{Tezduyar1991}, could then be reinterpreted as residual-based fine-scale models \cite{Brezzi1997b,Donea2003,Hughes2004b}. 
In addition, Hughes and collaborators established the VMS method as a framework for large-eddy simulations \cite{Hughes1999}. 
While initially fine-scale models used in this framework were based on eddy viscosity assumptions \cite{Hughes2003,Hughes2001a,Hughes2001,Munts}, more recently residual-based representations have been used for the subgrid scales \cite{Bazilevs2006,Bazilevs2007,Calo2004}. In the wider context of finite element discretizations of the Navier-Stokes equation, the VMS method in its residual-based form has also frequently been employed as a stabilization technique \cite{Hsu2012a,Hsu2014,Kamensky2015,Xu2016}.

Discontinuous Galerkin (DG) methods, first proposed by Reed and Hill in 1973 \cite{Reed1973}, have established themselves as an important paradigm for higher-order flow analysis over the last two decades \cite{Mavriplis:09.1,Wang:13.2}. The success of DG methods is based on a number of attractive properties, such as a rigorous mathematical foundation, the ability to use arbitrary orders of basis functions on general unstructured meshes, and a natural stability property for advective operators \cite{Cockburn:00.1,Hartmann:08.1}. Furthermore, DG methods are locally conservative, offer straightforward \textit{hp}-adaptivity, and are well-suited for parallel computing. For diffusion operators, there exist a variety of numerical flux formulations to tie discontinuous elements together, e.g., symmetric and non-symmetric interior penalty methods \cite{Arnold2001,Schillinger2016}, the local DG method \cite{Cockburn:05.1,Xu2009}, and the compact DG method \cite{Peraire2007}. A significant disadvantage of DG methods is the proliferation and increased global coupling of degrees of freedom with respect to standard continuous Galerkin methods. A remedy is the concept of hybridization \cite{Cockburn2009,Huerta2013,Kirby2012,Lehrenfeld:15.1}, where additional unknowns on element interfaces are introduced that reduce global coupling and facilitate static condensation strategies. 



In this paper, we explore a new residual-based variational multiscale formulation in the context of discontinuous Galerkin methods. In the past, several authors have investigated the use of the VMS method in a DG context. 
Examples are the multiscale discontinuous Galerkin methods introduced in \cite{Bochev2005,Buffa2006,Hughes2006} and methods for constructing discontinuous fine-scale bubble functions \cite{Coley2017,Sangalli2004}. These methods, however, maintain a continuous solution space for the coarse-scale problem and only use discontinuous subgrid scales for the fine-scale problem to capture element boundary layer effects. 
This approach is thus fundamentally different from the original VMS idea that we follow in this paper, that is, the decomposition of the true solution into a discontinuous coarse-scale function space and an accompanying discontinuous fine-scale function space. 
While several authors have investigated the enhancement of DG methods with fine-scale eddy viscosity models \cite{Collis2002,Collis2015,Ramakrishnan2004a}, DG methods based on a residual-based VMS subgrid-scale model are still largely unexplored. To some extent, this may be attributed to the importance of coarse-scale continuity
in the derivation of the VMS method, highlighted for example in \cite{Akkerman2008,Hughes1998,Hughes2004b}. 
The methodology that we develop in this paper no longer relies on the level of continuity of the coarse-scale function space. 



Our paper is structured as follows: in \Cref{sec:2}, we derive a general variational multiscale formulation that can accommodate discontinuous basis functions, including element interior (``volumetric'') and element boundary (``interface'') components of the fine scales. 
In \Cref{sec:3}, we transfer the fine-scale volumetric component into a residual-based model suitable for DG discretizations. 
In \Cref{sec:4}, we consider the one-dimensional Poisson problem to experimentally verify our new formulation for DG discretizations with linear finite elements, highlighting the importance of the fine-scale contributions. We demonstrate that existing discontinuous Galerkin formulations can be naturally rederived by particular fine-scale models in our VMS formulation. 
In this context, we discuss the fine-scale models associated with the interior penalty method. We also provide initial thoughts and some numerical experiments related to the extension of our formulation to multi-dimensional problems and higher-order DG discretizations. In \Cref{sec:5}, we consider the one-dimensional steady advection-diffusion problem, 
where we use our VMS formulation to investigate the relationship between upwinding and volumetric fine-scale models. We show that upwind numerical fluxes can be interpreted as an ad hoc remedy for missing volumetric fine-scale terms.






\section{Variational multiscale formulation in a discontinuous approximation space}
\label{sec:2}

In this section, we derive a general variational multiscale formulation that is able to accommodate discontinuous function spaces. In the scope of this article, we restrict ourselves to linear second-order PDEs. 
A boundary value problem with such a PDE is:
\begin{align}
\begin{cases}
\begin{alignedat}{3}
\mathcal{L} u &=f \qquad &&\text{in } \,\Omega \subset \mathbb{R}^d\\
u &= u_D \qquad&&\text{on }  \partial\Omega
\end{alignedat}\label{defPDE}
\end{cases}
\end{align}
with domain $\Omega$, dimension $d$, domain boundary $\partial\Omega$, source term $f$, and Dirichlet boundary data $u_D$. The differential operator $\mathcal{L} $ acts on the scalar solution $u$. We note that an extension of the following derivations to Neumann boundary conditions is straightforward. We assume throughout this paper that the exact solution of \eqref{defPDE} has continuous first derivatives, that is:
\begin{align}
u\in \mathcal{C}^1
\label{uinC1}
\end{align}

Our goal is to use a discontinuous Galerkin method to find a finite element solution to problem \eqref{defPDE}. To allow the use of discontinuous approximation functions later~on, we enlarge the space with respect to \eqref{uinC1} and consider solutions 
in the space of per-element continuous functions:
\begin{align}
\mathcal{V}(g) = \{ v \in L^2(\Omega) : v\big|_{K}\in \mathcal{C}^0 \,\, \forall\,K\in\mathcal{T} ,\,\, v=g \text{ on } \partial\Omega \}
\label{V}
\end{align}
where $K$ indicates an element, and $\mathcal{T}=\{K\}$ is the computational mesh.


\begin{table}[b!]
\centering
\begin{tabular}{lcl}\hline\hline\\[-0.2cm]
Jump operator               	& $\jump{w}\quad$        & $=w^{+} n^{+}+w^{-} n^{-}$                                          \\[0.3cm]
Average operator            	& $\avg{w}\quad$           & $=\frac{1}{2} \left( w^{+}+w^{-} \right)$                                     \\[0.3cm]
Volume $L^2$-inner product              	& $\big(w,u\big)_{K}\,\,$ & $=\int\limits_K w\, u $                                                    \\[0.5cm]
Surface $L^2$-inner product              	& $\big\langle w\, n,u\big\rangle_{\partial K}$ & $=\int\limits_{\partial K} w\, n \,u $                                                    \\[0.5cm]
Numerical interior domain   	& $\Omega_K \,$          &  $=\raisebox{2pt}{\scalebox{0.9}{$\bigcup\limits_{K\in\mathcal{T}}$}} \, K $\\[0.4cm]
                                 &                &$\big(w,u\big)_{\Omega_K} = \sum\limits_{K\in\mathcal{T}} \big(w,u\big)_{K}$ \\[0.5cm]
\multirow{1}{*}{\begin{tabular}[c]{@{}l@{}}Interior facets\\ \hspace{0.3cm}(Excludes domain boundary)\end{tabular}}  & $\Gamma \,\, $ & $=\raisebox{-1pt}{$\Big[$}\raisebox{2pt}{\scalebox{0.9}{$\bigcup\limits_{K\in\mathcal{T}}$}} \,\partial K \raisebox{-1pt}{$\Big]$}\!\setminus\! \partial \Omega$\\[0.4cm]
                                 &                &$\big\langle\jump{w},1\big\rangle_{\Gamma}\,\, = \sum\limits_{K\in\mathcal{T}} \big\langle w\, n,1\big\rangle_{\partial K \setminus \partial \Omega}$ \\[0.4cm]
                                 &                &$\big\langle\avg{w},1\big\rangle_{\Gamma} = \sum\limits_{K\in\mathcal{T}} \big\langle \frac{1}{2}w,1\big\rangle_{\partial K \setminus \partial \Omega}$ \\[0.4cm] \hline\hline\\[-0.3cm]
\end{tabular}
\caption{Collection of frequently used definitions.}
\label{tab:defs1}
\end{table}

The weak formulation of \cref{defPDE} is obtained by means of the method of weighted residuals, with Dirichlet constraints being enforced strongly. Using the notation defined in \cref{tab:defs1}, the weak formulation reads as follows:
\begin{align}
\begin{split}
&\text{Find }u\in \mathcal{V}(u_D) \text{ s.t.:}\\
&\quad \begin{cases}
\sum\limits_{K\in\mathcal{T}} \big( w, \mathcal{L} u   \big)_{K} = \sum\limits_{K\in\mathcal{T}} \big( w, f \big)_{K} \qquad \forall\, w \in \mathcal{V}(0) \\[0.4cm]
\begin{alignedat}{3}
&\jump{u} &&= 0&&\hspace{3.85cm} \text{on } \raisebox{2pt}{\scalebox{0.9}{$\bigcup\limits_{K\in\mathcal{T}}$}} \,\partial K \!\setminus\! \partial \Omega\\[0.2cm]
&\jump{\nabla u} &&= 0&&\hspace{3.85cm} \text{on } \raisebox{2pt}{\scalebox{0.9}{$\bigcup\limits_{K\in\mathcal{T}}$}} \,\partial K \!\setminus\! \partial \Omega
\end{alignedat}
\end{cases}
\label{uweakform}
\end{split}
\end{align}
where $w$ is a test function such that $w=0$ on $\partial \Omega$. 

The discontinuous nature of $\mathcal{V}(\cdot)$ requires that the weak formulation is defined on a per-element basis. This ensures the existence of derivatives, which would otherwise lead to Dirac layers at element boundaries in the upcoming derivations. Transmission conditions in the second and third line of \eqref{uweakform} couple the solution from element to element, thereby ensuring the continuity and uniqueness of the solution. 


The variational multiscale method suggests the split of the solution into a coarse-scale (finite element) component and a complementary fine-scale component. The following is a typical discontinuous Galerkin function space, which serves as the basis for the coarse-scale solution:
\begin{align}
\widebar{\mathcal{V}}(g) = \{ v \in L^2(\Omega) : v\big|_{K}\in P^p  (K)\,\, \forall\,K\in\mathcal{T} ,\,\, v=g \text{ on } \partial\Omega \}
\label{Vb}
\end{align}
$P^p(K)$ denotes the space of polynomial functions, defined on element~$K$, up to order~$p$. Notice that the function space $\widebar{\mathcal{V}}(g)$ satisfies
\begin{align}
\widebar{\mathcal{V}}(g) \subset  \mathcal{V}(g)
\end{align}

The goal of the variational multiscale method is to find a coarse-scale solution~$\bar{u}$, defined as some projection of  the true solution $u$ onto the space $\widebar{\mathcal{V}}(\cdot)$. The difference between the true solution $u$ and the coarse-scale solution $\bar{u}$ is defined as the fine-scale solution $u'$. We denote the projector used to define $\bar{u}$ as $\mathcal{P}$, so $\mathcal{P}:\mathcal{V}(\cdot) \rightarrow \widebar{\mathcal{V}}(\cdot)$. This projector is assumed to be a linear mapping. By definition, the projection $\mathcal{P}$ is also idempotent, that is: $\mathcal{P}(\mathcal{P}u) = \mathcal{P}u$. We obtain the following definitions:
\begin{align}
& \bar{u} \,\equiv \mathcal{P}u \in \widebar{\mathcal{V}}(\cdot) \label{defub}\\
& u' \equiv u-\bar{u} \,\,\Rightarrow\,\, u =\bar{u}+u'\label{defup}\\
& \mathcal{P} u' =  \mathcal{P} (u-\bar{u}) =  \mathcal{P} u - \mathcal{P}(\mathcal{P}u) = 0 \,\,\Rightarrow\,\, u' \in \ker(\mathcal{P}) \subset \mathcal{V}(\cdot) \label{upkerP}
\end{align}
Note that $u'$ is an element of the space $\ker(\mathcal{P})$, which we therefore denote the fine-scale space: 
\begin{align}
\mathcal{V}'(g) = \{ v \in V(g) : \mathcal{P}v = 0 \}
\label{Vp}
\end{align}
Functions in $\mathcal{V}'(\cdot)$ may be discontinuous across element boundaries.

The final assumption on $\mathcal{P}$ is that it is a surjective mapping. This means that for every $\bar{v}\in \widebar{\mathcal{V}}(\cdot)$, there exists at least one $v\in\mathcal{V}(\cdot)$ such that $\mathcal{P}v = \bar{v}$. Since $\widebar{\mathcal{V}}\subset\mathcal{V}$, and $\mathcal{P}$ is idempotent, this implies that $\mathcal{P}\bar{v} = \bar{v}$ for every $\bar{v}\in \widebar{\mathcal{V}}(\cdot)$. By construction of $\mathcal{V}'$, this means that $\widebar{\mathcal{V}}\cap\mathcal{V}' = \{0\} $. As a consequence, the spaces $\widebar{\mathcal{V}}$ and $\mathcal{V}'$ form a direct sum decomposition of $\mathcal{V}$:
\begin{align}
\mathcal{V} = \widebar{\mathcal{V}}\oplus \mathcal{V}'
\label{directSum}
\end{align}
Due to \eqref{directSum}, any true solution $u\in\mathcal{V}$ maps \textit{uniquely} into a coarse-scale solution $\bar{u}\in\widebar{\mathcal{V}}$ and a fine-scale solution $u'\in\mathcal{V}'$. This property is important for the well-posedness of the variational multiscale formulation.

With these definitions at hand, the variational multiscale approach can be used to obtain the variational discontinuous Galerkin formulation. First, the transmission conditions are written in terms of coarse-scale and fine-scale solutions:
\begin{align}
\begin{alignedat}{10}
&\jump{u} &&= 0 = \,\,\,\, \jump{\bar{u}+u'} &&= \,\,\,\, \jump{\bar{u}}+\jump{u'} && \quad \Rightarrow \quad \jump{\bar{u}} &&= -\jump{u'} &&\qquad \text{on } \Gamma\\
&\jump{\nabla u} &&= 0 = \jump{\nabla \bar{u}+\nabla u'} &&= \jump{\nabla \bar{u}}+ \jump{\nabla u'} &&\quad \Rightarrow\quad  \jump{\nabla\bar{u}} &&= -\jump{\nabla u'}&&\qquad \text{on } \Gamma
\label{ubup_cont}
\end{alignedat}
\end{align}
where $\Gamma$ denotes the interior facets (see \cref{tab:defs1}).

In the next step, definitions \cref{defub,defup,upkerP} and relations \cref{ubup_cont} are substituted into the weak form~\cref{uweakform}. We assume that Dirichlet boundary conditions can be perfectly represented in the coarse-scale function space. Therefore, the fine-scale solution equals zero on the domain boundary. We find the following variational formulation: 
\begin{align}
\begin{split}
&\text{Find }\bar{u},u'\in \widebar{\mathcal{V}}(u_D)\times \mathcal{V}'(0) \text{ s.t.:}\\
&\quad \begin{cases}
\big( \bar{w}\,, \mathcal{L} (\bar{u}+u')   \big)_{\Omega_K} \,= \big( \bar{w}\,, f \big)_{\Omega_K} \qquad \forall\, \bar{w} \in \widebar{\mathcal{V}}(0)\\
\big( w',  \mathcal{L} (\bar{u}+u') \big)_{\Omega_K} = \big( w' , f \big)_{\Omega_K} \qquad \forall\, w' \in \mathcal{V}(0)\\[0.1cm]
\begin{alignedat}{3}
&\jump{\bar{u}} &&= -\jump{u'}       \hspace{2.5cm}  &&\qquad \text{on } \Gamma\\[0.1cm]
&\jump{\nabla\bar{u}} &&= -\jump{\nabla u'}  &&\qquad \text{on } \Gamma
\end{alignedat}
\end{cases}
\label{ubupweakform}
\end{split}
\end{align}
where $\Omega_K$ denotes the set of element interior domains (see \cref{tab:defs1}).

The first line in \cref{ubupweakform} represents the variational coarse-scale formulation, which will be the basis for the VMS finite element discretization. It includes the effect of the fine scales on the coarse-scale solution, which will have to be modeled to close the formulation. The second line represents the fine-scale problem, which will be the basis of the fine-scale volumetric model. We deliberately keep the space of fine-scale test functions broader than the space of fine-scale solution functions. This will help us in the development of fine-scale models in \Cref{sec:3}.

We can then rewrite the variational coarse-scale formulation as follows:
\begin{align}
\begin{split}
&\text{Find }\bar{u}\in \widebar{\mathcal{V}}(u_D) \text{ s.t.:}\\
&\quad B(\bar{w},\bar{u}) \,+ s(\bar{w},\bar{u} ;\Gamma) \,+ \big(\mathcal{L}^*\! \bar{w}\,,  u' \big)_{\Omega_K}  \,+ k(\bar{w},u' ;\Gamma) \,= L(\bar{w}) \,\,\quad \forall \, \bar{w} \,\in \widebar{\mathcal{V}}(0)
\end{split}
\label{ubweakform}
\end{align}
where $B(\,\cdot\,,\,\cdot\,)$ and $L(\,\cdot\,)$ represent the bilinear and linear forms typically found in finite element formulations of the PDE at hand. 
The differential operator $\mathcal{L}^*$ is the adjoint of $\mathcal{L}$. The term $s(\,\cdot\,,\,\cdot\,;\Gamma)$  contains surface terms on $\Gamma$ that we obtain by performing integration by parts 
on $(\bar{w},\mathcal{L}\bar{u})_{\Omega_K}$. Similarly, the term $k(\,\cdot\,,\,\cdot\,;\Gamma)$ appears when we transform $(\bar{w},\mathcal{L}u')_{\Omega_K}$ into $(\mathcal{L}^*\! \bar{w}, u')_{\Omega_K}$. 

The variational coarse-scale formulation \cref{ubweakform} includes sums over all elements in the mesh, but the solution on each element does not yet communicate with the solution in other elements. 
To achieve element coupling, the fine-scale solution $u'$ on an element interface is re-expressed as follows:
\begin{align}
\begin{alignedat}{6}
&u'^{\pm}\, n^\pm &&= \quad\frac{1}{2}(u'^+ \!\!+\! u'^- )\, n^\pm &&\!+\!\quad\frac{1}{2}(u'^+ n^+ \!\!+\! u'^- n^-) &&= \avg{u'}\, n^\pm\!+\!\frac{1}{2}\jump{u'} \\
&\nabla u'^{\pm}\!\cdot n^\pm &&= \frac{1}{2}(\nabla u'^+ \!\!+\! \nabla u'^- )\cdot n^\pm &&\!+\!\frac{1}{2}(\nabla u'^+ \!\!\cdot n^+ \!\!+\! \nabla u'^- \!\!\cdot  n^-) &&= \avg{ \nabla u'}\cdot n^\pm\!+\!\frac{1}{2}\jump{ \nabla u'}
\end{alignedat}\label{kdjfd}
\end{align}
where quantities with $+$ or $-$ refer to the edge of the current element and the edge of the neighboring element, respectively. The transmission conditions in the third and fourth lines of \cref{ubupweakform} couple the solution between elements. Substituting these conditions into \cref{kdjfd}, we can remove a number of dependencies on the fine-scale solution to find the following new relations:
\begin{align}
\quad\begin{alignedat}{6}
&u'^{\pm}\, n^\pm &&= \avg{u'}\, n^\pm-\frac{1}{2}\jump{\bar{u}} \\
&\nabla u'^{\pm}\cdot n^\pm &&= \avg{ \nabla u'}\cdot n^\pm-\frac{1}{2}\jump{ \nabla \bar{u}}
\end{alignedat}\label{step4Key}
\end{align}
When we substitute these relations into $k(\bar{w},u' ;\Gamma)$, the variational coarse-scale formulation \cref{ubweakform} becomes globally coupled. 

We note that at this stage, no simplifications or approximations have been introduced. Therefore, when solving \cref{ubweakform}, we obtain the exact coarse-scale solution $\bar{u}$ when the correct fine-scale contributions are incorporated. Recall that the split of the solution into coarse-scale and fine-scale components is defined by the projector~$\mathcal{P}$. Therefore, this projector is explicitly and inextricably tied to the remaining fine-scale terms. These terms can be divided into the following two components: the volumetric (intra-element) contribution represented by $(\mathcal{L}^*\! \bar{w}\,,  u' )_{\Omega_K}$ and the fine-scale interface (inter-element) contributions that appear by substitution of \cref{step4Key} into \cref{ubweakform}.

The fine-scale interface terms originate from the lack of continuity in the solution. As such, they are a fundamental part of the variational multiscale method in a discontinuous Galerkin framework. In a numerical implementation, these terms may be treated explicitly or implicitly. Explicit treatment suggests the substitution of predefined expressions of the form:
\begin{align}
&\begin{alignedat}{3}
&\avg{u'} && = \raisebox{-2pt}{\scalebox{1.2}{$\Phi$}}^{\text{\tiny E}} \\
&\avg{\nabla u'} &&= \raisebox{-1pt}{\scalebox{1.1}{$\Theta$}}^{\text{\tiny E}} 
\end{alignedat}\label{explFSSM}
\end{align}
We will use this approach to verify our VMS formulation in the next section. 
In practice, however, the average fine-scale quantities are unknowns, which suggests an implicit treatment. This implies that the fine-scale interface terms need to be written in terms of coarse-scale terms to obtain expressions of the form: 
\begin{align}
&\begin{alignedat}{3}
&\avg{u'} &&= \raisebox{-2pt}{\scalebox{1.2}{$\Phi$}}^{\text{\tiny I}} (\avg{\bar{u}},\jump{\bar{u}},\avg{\nabla\bar{u}},\jump{\nabla\bar{u}},\cdots)\\
&\avg{\nabla u'} &&= \raisebox{-1pt}{\scalebox{1.1}{$\Theta$}}^{\text{\tiny I}} (\avg{\bar{u}},\jump{\bar{u}},\avg{\nabla\bar{u}},\jump{\nabla\bar{u}},\cdots)
\end{alignedat} \label{implFSSM}
\end{align}

Much of the remaining body of this article will be dedicated to understanding the nature of the projectors that are enforced by implicit fine-scale interface models of the form~\cref{implFSSM}.



\section{Discontinuous residual-based modeling of the volumetric fine-scale term}
\label{sec:3}

In $\mathcal{C}^1$-continuous VMS methods, the complete scale interaction is captured by the volumetric fine-scale term. A key assumption in the classical continuous VMS model is that the fine-scale solution vanishes at element boundaries \cite{Akkerman2008,Hughes2004b}. For a discontinuous Galerkin formulation, this assumption is particularly unsuitable. The fine-scale solution at element boundaries is naturally tied to the jump of the coarse-scale solution through \cref{step4Key}. An appropriate volumetric fine-scale model must take the fine-scale discontinuities at element boundaries into account. In this section, we derive a solution of the fine-scale problem stated in the second line of \cref{ubupweakform} that incorporates the related nonhomogeneous element-boundary conditions.

The fine-scale weak formulation of \cref{ubupweakform} involves a sum over all elements. 
For a single element $K$, the variational fine-scale formulation is written in a form that introduces the residual of the coarse-scale solution as a forcing term. This results in: 
\begin{align}
\begin{split}
&\text{Find } u'\in \mathcal{V}'(0) \text{ s.t.:}\\
&\quad \big(w',\mathcal{L} u'\big)_{K} = \big(w',f\big)_{K}-\big(w',\mathcal{L} \bar{u}\big)_{K} = \big(w',\mathcal{R}_{\bar{u}} \big)_{K} \qquad \forall\, w' \in \mathcal{V}(0)  \label{resBased}
\end{split}
\end{align}
We can use Green's identities to rewrite the leftmost term of \eqref{resBased}, which results in:
\begin{align}
\begin{split}
&\text{Find } u'\in \mathcal{V}'(0) \text{ s.t.:}\\
&\quad \big(\mathcal{L}^*\! w', u'\big)_{K} + k(w',u';\partial K) = \big(w',\mathcal{R}_{\bar{u}} \big)_{K} \qquad \forall\, w' \in \mathcal{V}(0)  \label{GF1}
\end{split}
\end{align}
The resulting interface terms are summarized in $k(\,\cdot\,,\,\cdot\,;\partial K)$, which corresponds to the $k(\,\cdot\,,\,\cdot\,;\Gamma)$ in \cref{ubweakform}, although now acting on the element boundary $\partial K$.

The fine-scale volumetric model can be obtained by making use of the Green's function, defined as follows:
\begin{align}
\begin{cases}
\, g(x,y) \in \mathcal{V}(0)\\
\, \mathcal{L}^*\! g(x,y) = \delta_x \quad&\text{for } y\in K\\
\, g(x,y) = 0  \quad \qquad &\text{for } y\in \partial K
\end{cases}\label{gDef}
\end{align}

We note that the specific form of the Green's function depends on the PDE. The first line in \eqref{gDef} ensures that the Green's function is a suitable test function and may be substituted in \cref{GF1} in place of $w'$. The following identities hold, where the parameter of integration and differentiation is $y$: 
\begin{align}
\int\limits_{K} \mathcal{L}^*\! g(x,y)\, u' \,\text{d}y = \int\limits_{K} \delta_x \, u' \,\text{d}y =  u' =  - k(g(x,y),u';\partial K) + \int\limits_{K} g(x,y)\,\mathcal{R}_{\bar{u}}  \,\text{d}y \label{GF2}
\end{align}

The term $k(g(x,y),u';\partial K)$ involves integration on element boundaries, enforcing per-element fine-scale boundary conditions. The fine-scale solution $u'$ assumes boundary values that are prescribed by an explicit or implicit interface model \cref{explFSSM} or \cref{implFSSM} substituted into \cref{step4Key}. For the Poisson equation, the term $k(g(x,y),u';\partial K)$ reads: 
\begin{align}
k(g(x,y),u';\partial K) = \int_{\partial K} \nabla g(x,y)\cdot n \, u' \text{d}y
\label{exk}
\end{align}

These expressions are defined locally on each element, a favorable property for implementation in a finite element setting. However, the evaluation of the integrals is computationally expensive. In the variational coarse-scale formulation \cref{ubweakform}, the fine-scale solution in each element contributes only in a weighted sense, which enables further simplifications. In the next section, we will consider numerical experiments on a one-dimensional domain that use DG discretizations with linear basis functions. In this case, the volumetric fine-scale term can be written as:
\begin{align}
\big(\mathcal{L}^*\!\bar{w},u' \big)_K = \Big(\mathcal{L}^*\!\bar{w},- k( g(x,y),u';x_j ) - k( g'(x,y),u';x_{j+1} ) + \int\limits_{K} g(x,y)\, \mathcal{R}_{\bar{u}} \,\text{d}y  \Big)_{K} 
\end{align}
Here, $x_j$ and $x_{j+1}$ refer to the left and right node of a one-dimensional element. 

For linear basis functions and constant model parameters, the coarse-scale residual $\mathcal{R}_{\bar{u}}$ is constant on each element, $\mathcal{L}^*\!\bar{w}$ is constant, and the boundary values $u'(x_j)$ and $u'(x_{j+1})$ do not vary in space. These terms can therefore be taken out of their respective integrals and the integrals are reduced to multiplications with average Green's function quantities:
\begin{align}
\big(\mathcal{L}^*\!\bar{w},u' \big)_K = \big(\mathcal{L}^*\!\bar{w},\tau \mathcal{R}_{\bar{u}} + \gamma_0 u'(x_j) -  \gamma_1 u'(x_{j+1})  \big)_{K} \label{gResult}
\end{align}
where:
\begin{align}
\tau &= \frac{1}{|K|} \int\limits_{K} \int\limits_{K} g(x,y) \d{y}\dx \label{taudef}\\
\gamma &= \frac{1}{|K|} \int\limits_{K}\int\limits_{\partial K} \dd{y} g(x,y) \d{y}\dx  \label{gammadef}
\end{align}
The exact definition of $\gamma$ depends on the definition of the functional $k( g'(x,y),u';\partial K)$. Equation \cref{gammadef} is based on the example functional of equation \cref{exk}. For the one-dimensional case, the boundary integration in \cref{gammadef} may be written as the sum of two nodal values, hence the use of $\gamma_0$ and $\gamma_1$ in \cref{gResult}. The form shown in \cref{gammadef} corresponds to a multi-dimensional problem. It should be noted, however, that in this case $u'$ would not be constant across all of $\partial K$, and \cref{gResult} would always lead to an approximation.



\section{A VMS discontinuous Galerkin method for second-order elliptic problems}
\label{sec:4}

In this section, we employ the variational multiscale principles presented above for the solution of a boundary value problem based on a second order elliptic operator.  Following the notations in \eqref{defPDE}, we define a standard Poisson problem as follows: 
\begin{align}
\quad\begin{cases}
\begin{alignedat}{3}
- \Delta u &= f \qquad &&\text{in } \,\Omega\subset \mathbb{R}^d\\
u &= u_D \qquad &&\text{on } \partial\Omega \label{PoissonProb}
\end{alignedat}
\end{cases}
\end{align}

The Laplace operator is omnipresent in engineering, e.g., for modeling diffusion in transport problems or balance of momentum in elasticity. We note that its numerical treatment in a DG context had posed a challenge for many years until a number of breakthrough developments were achieved in the late 1990s (see e.g.\  \cite{Arnold2000,Arnold2001,Cockburn1999}).


\subsection{Variational multiscale formulation}
\label{ssec:4.1}

Based on the strategy discussed in \Cref{sec:2}, we define the following variational formulation on a per-element basis: 
\begin{align}
\begin{split}
&\text{Find } u\in \mathcal{V}(u_D) \text{ s.t.:}\\
&\quad\begin{cases}
\big(w, - \Delta u\big)_{\Omega_K} = \big(w,f\big)_{\Omega_K} \qquad \forall\, w\in\mathcal{V}(0) \hspace{3cm}\\[0.15cm]
\begin{alignedat}{3}
&\jump{u} &&= 0       \hspace{2.34cm}  &&\qquad \text{on } \Gamma\\[0.1cm]
&\jump{\nabla u} &&= 0  &&\qquad \text{on } \Gamma
\end{alignedat}
\end{cases}\label{WF29}
\end{split}
\end{align}
To obtain \cref{WF29}, we multiply \cref{PoissonProb} with a test function $w$, integrate over each element domain and sum over all elements. The transmission conditions in line two and three ensure the uniqueness of the solution.

The solution and test functions are split into coarse-scale and fine-scale components. For now, we leave the projector $\mathcal{P}$ that defines the fine-scale space $\mathcal{V}'(\cdot)$ unspecified. Using the split and exploiting the linearity of the Laplace operator, we can rewrite \eqref{WF29} into the following variational formulation:  
\begin{alignat}{3}
\begin{split}
&\text{Find }\bar{u}\text{, }u' \in \widebar{\mathcal{V}}(u_D)\times \mathcal{V}'(0) \text{ s.t.:}\\
&\quad\begin{cases}
\big(\bar{w},- \Delta\bar{u}\big)_{\Omega_K} +\big(\bar{w},- \Delta u' \big)_{\Omega_K} = \big(\bar{w},f\big)_{\Omega_K} \hspace{1.5cm} \forall \,\bar{w}\in\widebar{\mathcal{V}}(0)\\
\big(w',- \Delta u' \big)_{\Omega_K} = \big(w',f+ \Delta\bar{u}\big)_{\Omega_K} = \big(w',\mathcal{R}_{\bar{u}}\big)_{\Omega_K} \hspace{0.54cm} \forall \,w \in\mathcal{V}(0) \\[0.15cm]
\begin{alignedat}{3}
&\jump{\bar{u}} &&= -\jump{u'}       \hspace{4.9cm}  &&\qquad \text{on } \Gamma\\[0.1cm]
&\jump{\nabla\bar{u}} &&= -\jump{\nabla u'}  &&\qquad \text{on } \Gamma
\end{alignedat}
\label{split123}
\end{cases}
\end{split}
\end{alignat}

In order to provide a mechanism to introduce a fine-scale model later~on, we separate the fine-scale solution from the differential operator with the help of the following general form of Green's identity:
\begin{align}
\quad\big(\bar{w},\mathcal{L}u'\big)_{K} = \big(\mathcal{L}^*\! \bar{w},u'\big)_{K} + k(\bar{w},u';\partial K) \label{GreensIdentity}
\end{align}
For the Poisson problem, the identity $\mathcal{L}^* = \mathcal{L} = - \Delta$ holds, and  the term $k(\bar{w},u';\partial K)$ follows as:
\begin{align}
&k(\bar{w},u';\partial K)  =   -\big\langle \bar{w},\nabla u'\cdot n\big\rangle_{\partial K} + \big\langle \nabla \bar{w}\cdot n , u'\big\rangle_{\partial K}\\
&k(\bar{w},u';\Gamma)  =  \sum\limits_{K\in\mathcal{T}} -\big\langle \bar{w},\nabla u'\cdot n\big\rangle_{\partial K} + \big\langle \nabla \bar{w}\cdot n , u'\big\rangle_{\partial K}
\end{align}
Substituting \cref{GreensIdentity} into \cref{split123} and performing integration by parts on the first term, we find the variational coarse-scale formulation:
\begin{align}
\begin{split}
&\big(\nabla \bar{w},\nabla \bar{u}\big)_{\Omega_K} \,-\, \sum\limits_{K\in\mathcal{T}} \big\langle \bar{w},\nabla \bar{u}\cdot n\big\rangle_{\partial K} \\[-0.2cm]
&\hspace{3.3cm} -\big(\Delta \bar{w}, u' \big)_{\Omega_K} + k(\bar{w},u';\Gamma) = \big(\bar{w},f\big)_{\Omega_K}  \qquad \forall \,\bar{w}\in\widebar{\mathcal{V}}(0)
\label{coarsescale2}
\end{split}
\end{align}

Note that at this stage, the element boundary values of all variables, that is, $\bar{u}$, $\nabla\bar{u}$, $u'$, $\nabla u'$, are evaluated from their element interior definition. The PDE is ``contained'' in each element. Yet, inter-element coupling is required to solve the global system. As proposed in \Cref{sec:2}, this element coupling will be enforced through manipulation of the fine-scale boundary terms. To this end, we substitute the fine-scale boundary identities of \cref{step4Key} into the term $k(\,\cdot\,,\,\cdot\,;\Gamma)$ to find:
\begin{align}
\begin{alignedat}{2}
&k(\bar{w},u';\Gamma)  =  \sum\limits_{K\in\mathcal{T}}\Big[ -\big\langle \bar{w}\,n,\avg{\nabla u'}\big\rangle_{\partial K} +\big\langle \frac{1}{2}\bar{w},\jump{\nabla \bar{u}}\big\rangle_{\partial K} + \\[-0.2cm]
&\hspace{4.5cm}\big\langle \nabla \bar{w}\cdot n , \avg{u'}\big\rangle_{\partial K} - \big\langle \frac{1}{2}\nabla \bar{w}, \, \jump{\bar{u}} \big\rangle_{\partial K} \Big] \label{BCterms1}
\end{alignedat}
\end{align}
We note that all fine-scale contributions that appear in \cref{BCterms1} are single-valued on element interfaces. Therefore, their sum across all elements yields:
\begin{align}
\begin{split}
 k(\bar{w},u';\Gamma)  &=  - \big\langle \jump{\bar{w}},\avg{\nabla u'}\big\rangle_{\Gamma} +\big\langle \avg{\bar{w}},\jump{\nabla \bar{u}}\big\rangle_{\Gamma}  + \big\langle \jump{\nabla \bar{w}}, \, \avg{u'}\big\rangle_{\Gamma} - \big\langle \avg{\nabla \bar{w}} ,\jump{\bar{u}}\big\rangle_{\Gamma}   \label{BCterms2}
\end{split}
\end{align}
where we have collected neighboring $\bar{w}\, n$ and $\nabla\bar{w} \cdot n$ terms in \cref{BCterms1} as jumps and neighboring $\frac{1}{2}\bar{w}$ and $\frac{1}{2}\nabla \bar{w}$ terms as averages.

Substitution of \cref{BCterms2} into \cref{coarsescale2} results in the following globally coupled formulation: 
\begin{align}
\begin{split}
&\big(\nabla \bar{w},\nabla \bar{u}\big)_{\Omega_K} - \sum\limits_{K\in \mathcal{T}} \big\langle \bar{w},\nabla \bar{u}\cdot n\big\rangle_{\partial K} +\big\langle \avg{\bar{w}},\jump{\nabla \bar{u}}\big\rangle_{\Gamma}  - \big\langle \avg{\nabla \bar{w}} ,\jump{\bar{u}}\big\rangle_{\Gamma}  \\[-0.15cm]
&\hspace{0.4cm} -\big(\Delta \bar{w}, u' \big)_{\Omega_K} - \big\langle \jump{\bar{w}},\avg{\nabla u'}\big\rangle_{\Gamma} + \big\langle \jump{\nabla \bar{w}}, \, \avg{u'}\big\rangle_{\Gamma} = \big(\bar{w},f\big)_{\Omega_K}  \hspace{0.5cm}  \forall \,\bar{w}\in\widebar{\mathcal{V}}(0)
\label{form1}
\end{split}
\end{align}
The second and third terms in \eqref{form1} can be simplified according to:
\begin{align}
\small\begin{split}
&\big\langle\avg{\bar{w}},\jump{\nabla \bar{u}}\big\rangle_{\partial K} - \big\langle\bar{w}^{+}, \nabla \bar{u}^{+}\!\cdot n^{+}\big\rangle_{\partial K} -\big\langle\bar{w}^{-}, \nabla \bar{u}^{-}\!\cdot n^{-}\big\rangle_{\partial K}\\[0.1cm]
&\hspace{0.5cm}= \int\limits_{\partial K} \!\!\Big( \frac{1}{2}(\bar{w}^+ \!+ \bar{w}^-)(\nabla \bar{u}^+ \!\!\cdot\! n^+ \!+ \nabla \bar{u}^- \!\cdot n^-) - (\bar{w}\, \nabla \bar{u}\!\cdot n)^+ \!- (\bar{w}\,\nabla \bar{u}\!\cdot n)^- \Big) \\[-0.2cm]
&\hspace{1.5cm}= \int\limits_{\partial K} \!\frac{1}{2} \Big( - \bar{w}^+ n^+ \!\!\cdot\! \nabla \bar{u}^+  \!- \bar{w}^+  n^+ \!\!\cdot\! \nabla \bar{u}^- \!- \bar{w}^- n^- \!\!\cdot\! \nabla \bar{u}^+  \!- \bar{w}^- n^- \!\!\cdot\! \nabla \bar{u}^-   \Big)  \\[-0.2cm]
&\hspace{2.5cm}= - \int\limits_{\partial K} \!\frac{1}{2} ( \bar{w}^+ n^+  \!+ \bar{w}^- n^- ) ( \nabla \bar{u}^- \!+ \nabla \bar{u}^+ ) = -\big\langle\jump{\bar{w}},\avg{\nabla \bar{u}}\big\rangle_{\partial K}
\label{collection}
\end{split}
\end{align}
which yields the final variational coarse-scale formulation of the Poisson problem:
\begin{align}
\begin{split}
&\hspace{1.0cm}\big(\nabla \bar{w},\nabla \bar{u}\big)_{\Omega_K} \!- \big\langle\jump{\bar{w}},\avg{\nabla \bar{u}}\big\rangle_{\Gamma}  - \big\langle \avg{\nabla \bar{w}} ,\jump{\bar{u}}\big\rangle_{\Gamma}  \\
&\hspace{0cm} -\big(\Delta \bar{w}, u' \big)_{\Omega_K} \!- \big\langle \jump{\bar{w}},\avg{\nabla u'}\big\rangle_{\Gamma} + \big\langle \jump{\nabla \bar{w}}, \, \avg{u'}\big\rangle_{\Gamma} = \big(\bar{w},f\big)_{\Omega_K}  \qquad \forall \,\bar{w}\in\widebar{\mathcal{V}}(0) \label{coarsescale3}
\end{split}
\end{align}
Recall again that the projector $\mathcal{P}$ is formally required to construct $\mathcal{V}'(\cdot)$. The choice of this projector affects two parts of the obtained coarse-scale weak formulation \cref{coarsescale3}. These are the volumetric fine-scale term, and the two fine-scale element interface terms. Once these are prescribed by appropriate expressions, then the projector is implicitly defined, and the associated coarse-scale solution $\bar{u}=\mathcal{P}u$ will solve \cref{coarsescale3}.

It is interesting to note that  \eqref{coarsescale3} shares many similarities with classical DG formulations that have been developed in the last two decades. For illustration purposes, we show a number of classical DG methods for the Poisson problem in \Cref{tab:DGforms}. These are typically derived via a mixed method approach \cite{Arnold2001}. We observe that the first three terms in \eqref{coarsescale3} are exactly the ones proposed in the classical DG formulation by Bassi and Rebay in 1997 \cite{Bassi1997}. We also observe that they appear in almost every classical DG method summarized in \Cref{tab:DGforms}.  

\begin{table}[t]
\centering
\caption{Overview of classical DG formulations for the Poisson problem (adapted from \cite{Arnold2001}). We note that $r$, $l$ and $r_e$ denote lifting operators.}
\label{tab:DGforms}
\begin{tabular}{ll}
\hline\hline \\[-0.3cm]
\textbf{Method name }             & \textbf{Global weak formulation} \\[0.1cm] \hline\\[-0.1cm]
Bassi-Rebay \cite{Bassi1997}       & \small{$\big(\nabla \bar{w},\nabla \bar{u}\big)_{\Omega_K} \!\!\!-\! \big\langle\jump{\bar{w}},\avg{\nabla \bar{u}}\big\rangle_{\Gamma} \!-\! \big\langle\avg{\nabla \bar{w}} ,\jump{\bar{u}}\big\rangle_{\Gamma} \!+\! \big(r(\jump{\bar{w}}),r(\jump{\bar{u}})\big)_{\Omega_K}$} \\[0.3cm]
Brezzi \textit{et al.} \cite{Brezzi1997}& \specialcell[c]{\small{$\big(\nabla \bar{w},\nabla \bar{u}\big)_{\Omega_K} \!\!\!-\! \big\langle\jump{\bar{w}},\avg{\nabla \bar{u}}\big\rangle_{\Gamma} \!-\! \big\langle\avg{\nabla \bar{w}} ,\jump{\bar{u}}\big\rangle_{\Gamma} \!+\! \big(r(\jump{\bar{w}}),r(\jump{\bar{u}})\big)_{\Omega_K} $} \\[0.15cm] \small{$- \big\langle\jump{\bar{w}},\mu\,\avg{r_e (\bar{u})}\big\rangle_{\Gamma}$} } \\[0.6cm]
Local DG \cite{Cockburn1998}  & \specialcell[c]{\small{$\big(\nabla \bar{w},\nabla \bar{u}\big)_{\Omega_K} - \big\langle\jump{\bar{w}},\avg{\nabla \bar{u}}\big\rangle_{\Gamma} - \big\langle\avg{\nabla \bar{w}} ,\jump{\bar{u}}\big\rangle_{\Gamma} +\big\langle\jump{\bar{w}},\frac{\eta}{h}\, \jump{\bar{u}}\big\rangle_{\Gamma}$} \\[0.15cm] \small{$+\big\langle\beta\cdot \jump{\bar{w}},\jump{\nabla \bar{u}}\big\rangle_{\Gamma}+\big\langle\jump{\nabla \bar{w}},\beta\cdot \jump{\bar{u}}\big\rangle_{\Gamma}$} \\[0.15cm] \small{$+ \big(r(\jump{\bar{w}})+l(\beta\cdot\jump{\bar{w}}),r(\jump{\bar{u}})+l(\beta\cdot\jump{\bar{u}})\big)_{\Omega_K}$}} \\[0.9cm]
\small{Interior penalty}  \cite{Douglas1976}  & \small{$\big(\nabla \bar{w},\nabla \bar{u}\big)_{\Omega_K} - \big\langle\jump{\bar{w}},\avg{\nabla \bar{u}}\big\rangle_{\Gamma} - \big\langle\avg{\nabla \bar{w}} ,\jump{\bar{u}}\big\rangle_{\Gamma}+ \big\langle\jump{\bar{w}},\frac{\eta}{h}\, \jump{\bar{u}}\big\rangle_{\Gamma} $}\\[0.3cm]
Bassi \textit{et. al} \cite{Bassi1997b}    & \small{$\big(\nabla \bar{w},\nabla \bar{u}\big)_{\Omega_K} \!-\! \big\langle\jump{\bar{w}},\avg{\nabla \bar{u}}\big\rangle_{\Gamma} \!-\! \big\langle\avg{\nabla \bar{w}} ,\jump{\bar{u}}\big\rangle_{\Gamma} \!-\! \big\langle\jump{\bar{w}},\mu\,\avg{r_e (\bar{u})}\big\rangle_{\Gamma}$} \\[0.3cm]
{Baumann-Oden}  \cite{Baumann1999}   & \small{$\big(\nabla \bar{w},\nabla \bar{u}\big)_{\Omega_K} - \big\langle\jump{\bar{w}},\avg{\nabla \bar{u}}\big\rangle_{\Gamma} + \big\langle\avg{\nabla \bar{w}} ,\jump{\bar{u}}\big\rangle_{\Gamma}$} \\[0.3cm]
NIPG   \cite{Riviere1999}       & \small{$\big(\nabla \bar{w},\nabla \bar{u}\big)_{\Omega_K} - \big\langle\jump{\bar{w}},\avg{\nabla \bar{u}}\big\rangle_{\Gamma} + \big\langle\avg{\nabla \bar{w}} ,\jump{\bar{u}}\big\rangle_{\Gamma}+ \big\langle\jump{\bar{w}},\frac{\eta}{h}\, \jump{\bar{u}}\big\rangle_{\Gamma}$} \\[0.3cm]
{Babu\v{s}ka-Zl\'amal}  \cite{Babuska1973} & \small{$\big(\nabla \bar{w},\nabla \bar{u}\big)_{\Omega_K} + \big\langle\jump{\bar{w}},\frac{\eta}{h}\, \jump{\bar{u}}\big\rangle_{\Gamma}$} \\[0.3cm]
Brezzi \textit{et al.} \cite{Brezzi2000} & \small{$\big(\nabla \bar{w},\nabla \bar{u}\big)_{\Omega_K} - \big\langle\jump{\bar{w}},\mu\,\avg{r_e (\bar{u})}\big\rangle_{\Gamma}$} \\[0.3cm]
\hline\hline 
\end{tabular}
\end{table}

\subsection{Numerical experiments with linear basis functions in 1-D}
\label{ssec:4.2}
In this section, we investigate the effect of the fine-scale interface terms in \cref{coarsescale3}, that is, $\avg{u'}$ and $\avg{\nabla u'}$.  To this end, we remove all volumetric fine-scale terms by restricting ourselves to linear trial and test functions in \cref{coarsescale3}. Due to  $\mathcal{L}^*\!\bar{w}=0$ in this case, only the fine-scale interface terms remain, while all volumetric intra-element terms vanish. As a consequence, the projector $\mathcal{P}$ is wholly defined by (or, conversely, completely defines) the fine-scale interface terms. We show this experimentally by substituting explicit models in place of the fine-scale interface terms, see \cref{explFSSM}. 


When $\bar{u}$ is constructed using an $H^1$ projection, it is nodally coincident with $u$. To obtain a closed coarse-scale formulation, we use this condition to simplify the terms $\avg{u'}$ and $\avg{\nabla u'}$. One can directly conclude that \mbox{$\avg{u'}=0$} on element interfaces. With nodal exactness, an explicit formulation for $\avg{\nabla u'}$ can be obtained as:
\begin{align}
\avg{\nabla u'}\Big|_{\hat{x}}& = \nabla u(\hat{x})-\avg{\nabla \bar{u}}\Big|_{\hat{x}} = \dd{x} u(x)\Big|_{\hat{x}}-\frac{1}{2}\left(\frac{u({\hat{x}})-u({\hat{x}}-h)}{h}+\frac{u({\hat{x}}+h)-u({\hat{x}})}{h}\right) \nonumber\\
&= \dd{x} u(x)\Big|_{\hat{x}}-\frac{u({\hat{x}}+h)-u({\hat{x}}-h)}{2h} \label{nabu}
\end{align}
where ${\hat{x}}$ is a coordinate on the element boundary. We assume that the element width $h$ is equal on neighboring elements. 

\begin{figure}[ht]
\centering
\subfloat{ \begin{tabular}{c}\vspace{-6.3cm}\\
\textbf{Model problem} \\
$f=2$ \\
$x_0 = 0$, $x_1 = 5$\\
$u_0 = 1$, $u_1 = 3$ \\\\
\textbf{Exact solution} \\
$u(x) = x^2+\frac{26}{5}x+1$   \\ \\
\textbf{Fine-scale model}\\
For ${\hat{x}}\in \Gamma$:  \\[0.1cm]
$\avg{u'}\big|_{\hat{x}}=\raisebox{-2pt}{\scalebox{1.2}{$\Phi$}}^{\text{\tiny E}}\big|_{\hat{x}} =  0$\\[0.2cm]
$\avg{\nabla u'}\big|_{\hat{x}}= \raisebox{-1pt}{\scalebox{1.1}{$\Theta$}}^{\text{\tiny E}}\big|_{\hat{x}} = 0$
\end{tabular}
}\hspace{0.2cm}
\subfloat{ %
  \includegraphics[width=0.62\textwidth]{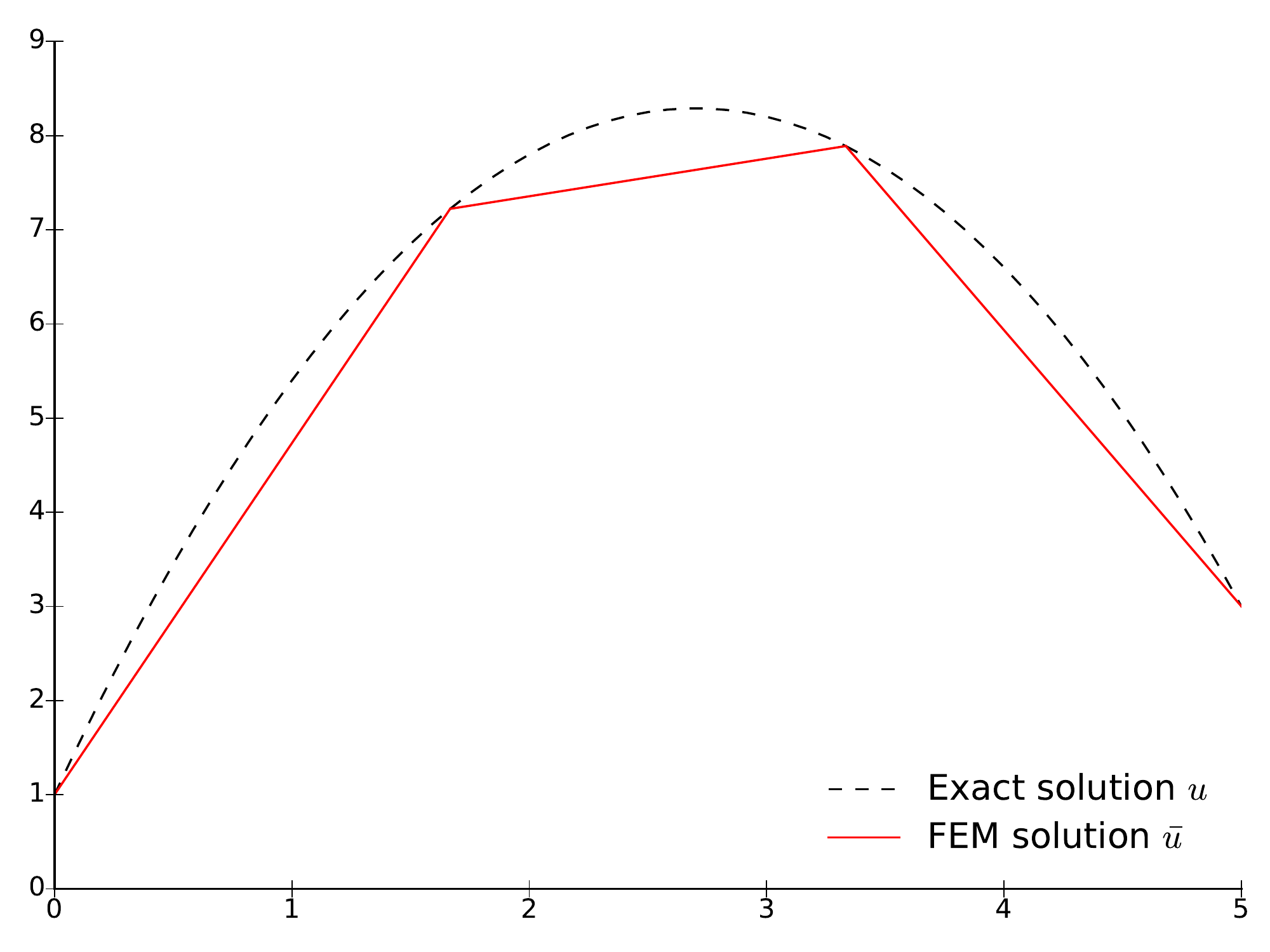}%
}
\caption{Numerical experiment 1: for linear basis functions and constant $f$, the coarse-scale solution $\bar{u}$ is nodally exact, when all fine-scale terms are omitted.} 
\label{fig:exp1}
\end{figure}

In a one-dimensional Poisson problem, a constant force $f$ yields a parabolic exact solution. For an arbitrary parabola, \mbox{$u=ax^2+bx+c$}, it follows that the last term in \cref{nabu} is equal to the derivative of the true solution, \mbox{$\dd{x}u = 2ax+b$}, according to:
\begin{align}
\frac{u(x+h)-u(x-h)}{2h}  &=  \frac{ax^2 +2axh+ah^2+bx+bh-ax^2 +2axh-ah^2-bx+bh}{2h} \nonumber\\ 
&= 2ax+b =\dd{x} u(x)  \nonumber
\end{align}
By substituting this relation into \cref{nabu}, we can show that for a constant force $f$ and linear basis functions, it holds that $\avg{\nabla u'}=0$. An inverse argument says that when all fine-scale terms are omitted from the variational coarse-scale formulation \eqref{coarsescale3}, then the $H^1$ projection of $u$ must be obtained in case of a constant $f$ and linear basis functions. We verify this argument with the first numerical experiment illustrated in \cref{fig:exp1}, where we use three equidistant linear discontinuous finite elements.

\begin{figure}[ht]
\hspace{-0.1cm}\subfloat{%
\begin{tabular}{c}\vspace{-5.8cm}\\
\textbf{Model problem} \\
$f=10(x-x^2)$ \\
$x_0 = 0$, $x_1 = 1$\\
$u_0 = 0$, $u_1 = 0.1$ \\\\
\textbf{Exact solution} \\
$u(x) = -\frac{5}{3}x^3+\frac{10}{12} x^4+\frac{14}{15}x$  \\ \\
\textbf{Fine-scale model}\\
For ${\hat{x}}\in \Gamma$:  \\
$\raisebox{-2pt}{\scalebox{1.2}{$\Phi$}}^{\text{\tiny E}}\big|_{\hat{x}}= 0$\\
$\raisebox{-1pt}{\scalebox{1.1}{$\Theta$}}^{\text{\tiny E}} \big|_{\hat{x}}=\nabla u({\hat{x}})-\avg{\nabla \bar{u}}\Big|_{\hat{x}}$\\[0.2cm]
$= \dd{x} u(x)\Big|_{\hat{x}}-\frac{u({\hat{x}}+h)-u({\hat{x}}-h)}{2h}$
\end{tabular}
}\hspace{-0.2cm}
\subfloat{
  \includegraphics[width=0.62\textwidth]{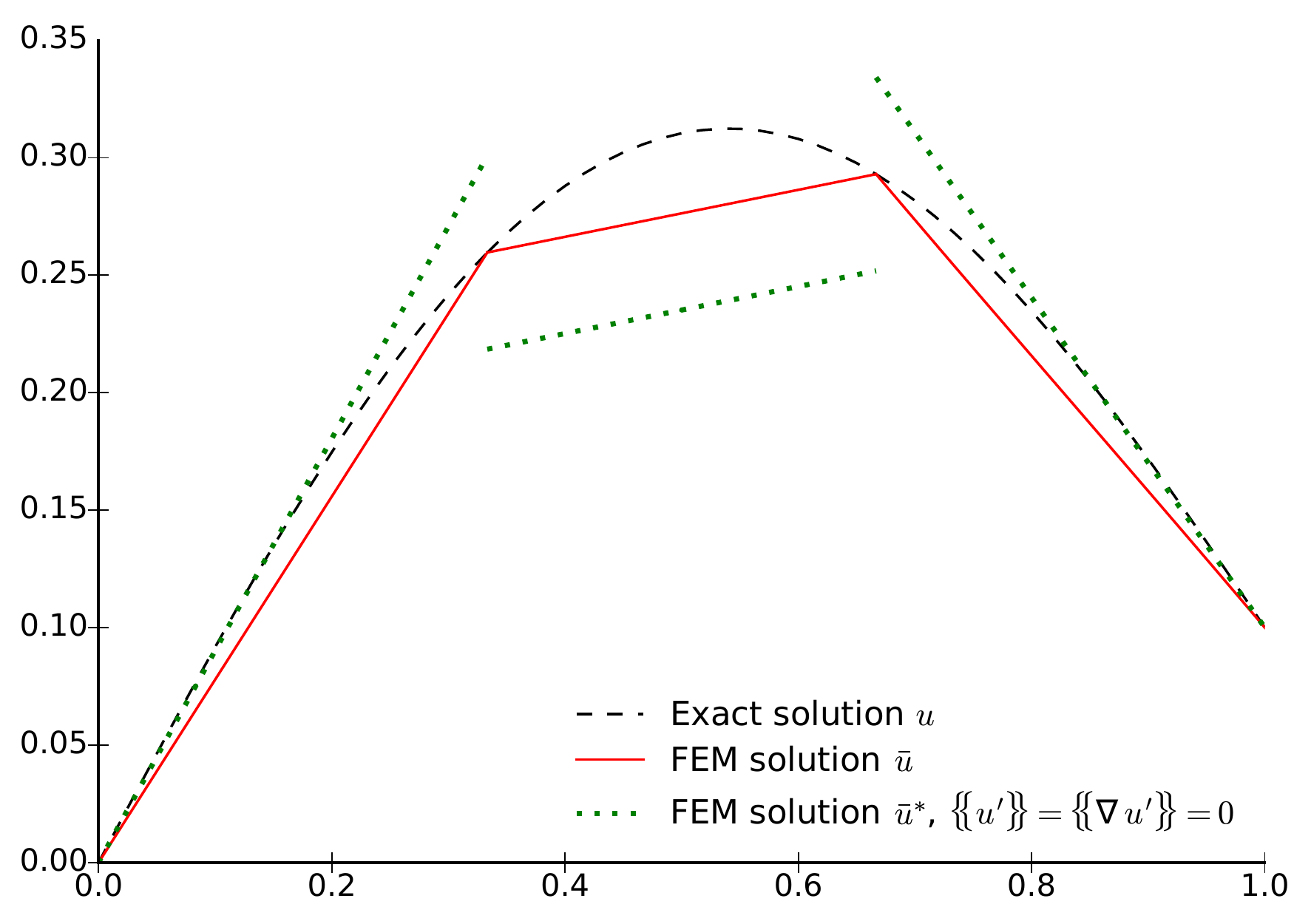}
}
\caption{Numerical experiment 2: for linear basis functions and a given parabolic $f$, the coarse-scale solution $\bar{u}$ is nodally exact when the correct $\avg{\nabla u'}$ is substituted as an explicit interface model.}%
  \label{fig:exp2} 
\end{figure}

For an arbitrary $f$, the nodally exact $H^1$ projection of $u$ does not necessarily imply $\avg{\nabla u'}=0$. In this case, nodal coincidence can still be achieved by using the explicit formulation for $\avg{\nabla u'}$ from \cref{nabu}, such that we obtain additional 
contributions to the force vector. In the second numerical experiment shown in \cref{fig:exp2}, we verify this argument for an example parabolic force, $f=10(x-x^2)$. In addition, we plot the result when the fine-scale terms are omitted, which yields a solution that is not nodally exact and exhibits significant jumps at element interfaces. 

\begin{figure}[ht]
\hspace{-0.2cm}\subfloat{
\begin{tabular}{c}\vspace{-6.1cm}\\
\textbf{Model problem} \\
$f=10(x-x^2)$ \\
$x_0 = 0$, $x_1 = 1.5$\\
$u_0 = 0$, $u_1 = 0.1$ \\\\
\textbf{Exact solution} \\
$u(x) = -\frac{5}{3}x^3+\frac{10}{12} x^4+\frac{241}{240}x$    \\ \\
\textbf{Fine-scale model}\\
For $\hat{x}\in \Gamma$:  \\
$\raisebox{-2pt}{\scalebox{1.2}{$\Phi$}}^{\text{\tiny E}}\big|_{\hat{x}}=u(\hat{x}) - \avg{\mathcal{P}_{\! L^2}u}\Big|_{\hat{x}}$\\
$\raisebox{-1pt}{\scalebox{1.1}{$\Theta$}}^{\text{\tiny E}} \big|_{\hat{x}}=\nabla u(\hat{x})-\avg{\nabla \,\mathcal{P}_{\! L^2}u }\Big|_{\hat{x}}$
\end{tabular}
}\hspace{-0.2cm}
\subfloat{%
  \includegraphics[width=0.6\textwidth,trim={0.2cm 0.cm 0.3cm 0.0cm},clip]{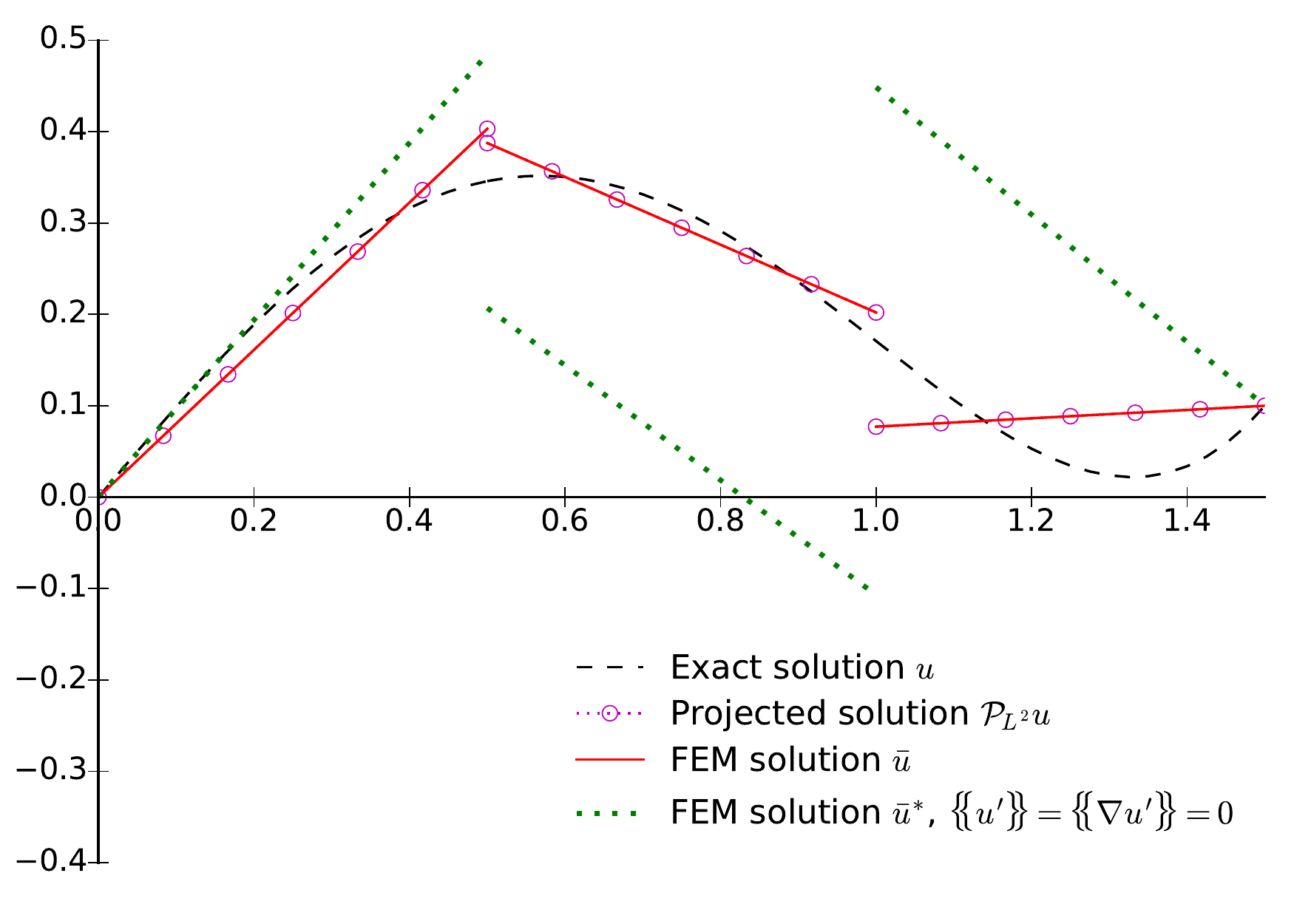}
}
  \caption{Numerical experiment 3: the $L^2$ projection of $u$ is retrieved when the correct $\avg{\nabla u'}$ and $\avg{u'}$ are substituted as explicit interface models.}%
  \label{fig:exp3}
\end{figure}

Nodal exactness is not necessarily a favorable quality for a discontinuous approximation space. The additional degrees of freedom carry no value. Often one is more interested in a solution that minimizes some error norm. For a given discretization, the $L^2$ error is minimized by the $L^2$ projection of $u$ onto the finite element space, that is, $\bar{u} = \mathcal{P}_{\! L^2} u$.  The third numerical experiment shown in \Cref{fig:exp3} verifies that explicit formulations for $\avg{u'}$ and $\avg{\nabla u'}$, which correspond to the desired $\bar{u}$, indeed lead to the $L^2$ projection of $u$. To find the reference solution, we solve a finite element problem to obtain $u_h = \mathcal{P}_{\! L^2} u$. From this finite element approximation, we extract the explicit formulation of the fine-scale interface terms, which become force contributions. We note that for our purpose, the $L^2$ projection must satisfy the Dirichlet boundary conditions, such that the fine-scale solution \mbox{$u'=u-\bar{u}\in \mathcal{V'}\big(0)$} is zero on the domain boundary. \Cref{fig:exp3} also includes the result of a DG discretization that omits all fine-scale terms in the variational coarse-scale formulation \eqref{coarsescale3}, which does not lead to an acceptable solution.

All three numerical experiments demonstrate that suitable fine-scale terms at element interfaces are essential for the optimal performance of DG discretizations. Explicit fine-scale interface models require prior knowledge of the exact solution, which is usually not available in practical applications. This is different for implicit interface models, which we will discuss next. 

\subsection{Multiscale interpretation of classical discontinuous Galerkin formulations}
\label{ssec:4.3}
In this section, we investigate the fine-scale interface models that naturally recover classical DG methods from the variational formulation~\eqref{coarsescale3}. 
Most of the classical DG formulations in \cref{tab:DGforms} can be written in the general form:
\begin{align}
\begin{split}
& \hspace{1cm} \big(\nabla \bar{w},\nabla \bar{u}\big)_{\Omega_K} \!- \big\langle\jump{\bar{w}},\avg{\nabla \bar{u}}\big\rangle_{\Gamma}  \!- \big\langle \avg{\nabla \bar{w}} ,\jump{\bar{u}}\big\rangle_{\Gamma} \,+ \\
& \big(\mathcal{F}(\jump{\bar{w}}), \mathcal{F}(\jump{\bar{u}}) \big)_{\Omega_K} \!\!-\! \big\langle \jump{\bar{w}},\raisebox{-2pt}{\scalebox{1.2}{$\Phi$}}^{\text{\tiny I}}(\bar{u})\big\rangle_{\Gamma} \!+\! \big\langle \jump{\nabla \bar{w}}, \raisebox{-1pt}{\scalebox{1.1}{$\Theta$}}^{\text{\tiny I}}(\bar{u}) \big\rangle_{\Gamma} = \big( \bar{w},f \big)_{\Omega_K}  \quad \forall \,\bar{w}\in\widebar{\mathcal{V}}(0)  
\end{split}\label{genClassform}
\end{align}
where $\mathcal{F}(\cdot)$ is linear and $\mathcal{F}(0) = 0$ \cite{Arnold2001}.
When we compare \cref{genClassform} to the variational coarse-scale formulation \cref{coarsescale3}, we observe that \cref{genClassform} is a special case of the multiscale formulation, where the fine-scale solution satisfies:
\begin{align}
\begin{split}
& -\big(\Delta \bar{w}, u' \big)_{\Omega_K} \!- \big\langle \jump{\bar{w}},\avg{\nabla u'}\big\rangle_{\Gamma} + \big\langle \jump{\nabla \bar{w}}, \, \avg{u'}\big\rangle_{\Gamma} =  \\
& \big(\mathcal{F}(\jump{\bar{w}}), \mathcal{F}(\jump{\bar{u}}) \big)_{\Omega_K} - \big\langle \jump{\bar{w}},\raisebox{-2pt}{\scalebox{1.2}{$\Phi$}}^{\text{\tiny I}}(\bar{u})\big\rangle_{\Gamma} + \big\langle \jump{\nabla \bar{w}}, \raisebox{-1pt}{\scalebox{1.1}{$\Theta$}}^{\text{\tiny I}}(\bar{u}) \big\rangle_{\Gamma}  \quad \forall \,\bar{w}\in\widebar{\mathcal{V}}(0)  \label{classformequal}
\end{split}
\end{align}
which may be expressed as:
\begin{align}
\begin{split}
&\hspace{1.1cm}-\big(\Delta \bar{w}, u'\big)_{\Omega_K} - \big(\mathcal{F}(\jump{\bar{w}}), \mathcal{F}(\jump{\bar{u}}) \big)_{\Omega_K}   \\
&  - \big\langle \jump{\bar{w}}, \avg{\nabla u'}-\raisebox{-2pt}{\scalebox{1.2}{$\Phi$}}^{\text{\tiny I}}(\bar{u}) \big\rangle_{\Gamma} + \big\langle\jump{\nabla \bar{w}}, \avg{u'}- \raisebox{-1pt}{\scalebox{1.1}{$\Theta$}}^{\text{\tiny I}}(\bar{u})\big\rangle_{\Gamma} =0\hspace{0.8cm} \forall \,\bar{w}\in\widebar{\mathcal{V}}(0)\label{appC1}
\end{split}
\end{align}

\begin{figure}[h!]
\centering
\subfloat[$\bar{w}$ corresponding to the $\jump{\nabla \bar{w}}$ terms.] {\includegraphics[width=0.48\textwidth]{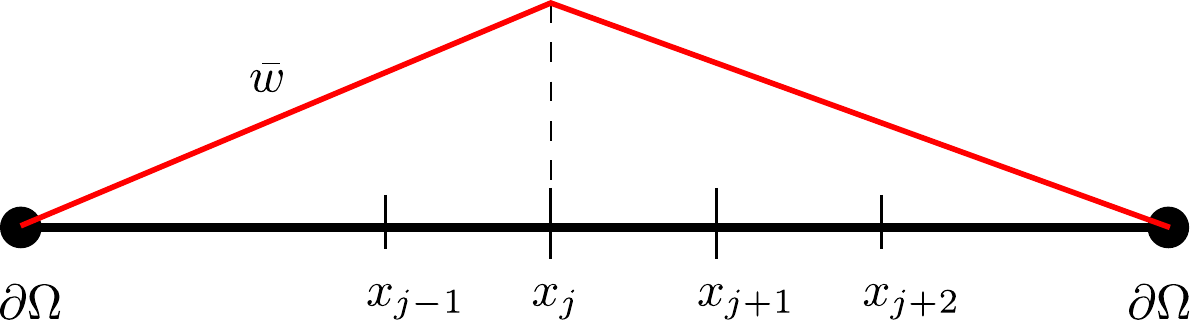} \label{fig:appC1}}  \hspace{0.1cm}
\subfloat[$\bar{w}$ corresponding to the $\jump{\bar{w}}$ terms.] {\includegraphics[width=0.48\textwidth]{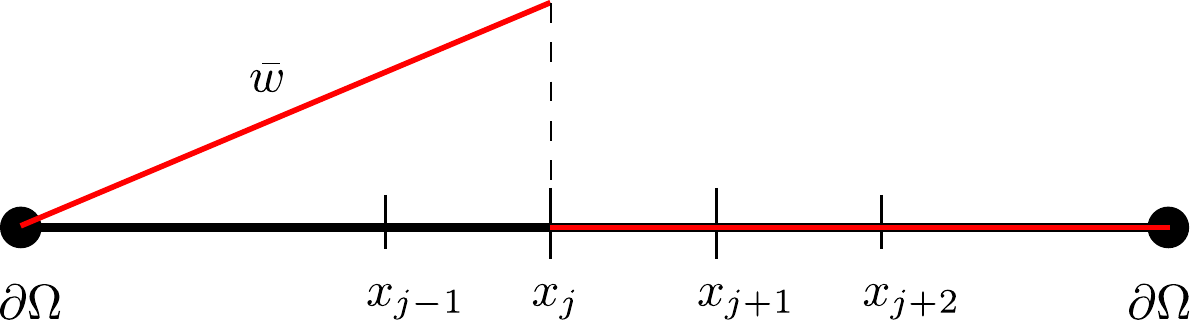} \label{fig:appC2}}
\caption{Choices of test functions to obtain the nodal identities.}
\end{figure}

For a one-dimensional problem, the coarse-scale test function $\bar{w}$ can be chosen as the continuous piecewise linear function shown in \cref{fig:appC1}. In this case, all but the last term in \cref{appC1} are zero. A similar test function can be constructed for each element interface. Therefore, it holds:
\begin{align}
\quad \big\langle \,1\, , \avg{u'} -\raisebox{-1pt}{\scalebox{1.1}{$\Theta$}}^{\text{\tiny I}}(\bar{u})  \big\rangle_{\hat{x}} = 0 \quad \Rightarrow \quad \avg{u'}\Big|_{{\hat{x}}}= \raisebox{-1pt}{\scalebox{1.1}{$\Theta$}}^{\text{\tiny I}}(\bar{u})\Big|_{{\hat{x}}}   &\qquad {\hat{x}}\in \Gamma \label{appC3}
\end{align}

If we choose the piecewise linear function shown in \cref{fig:appC2}, with a jump at some arbitrary element interface, as the test function $\bar{w}$, then the first term in \cref{appC1} can be removed. Since \cref{appC3} still holds, it follows that:
\begin{align}
\quad - \big(\mathcal{F}(1), \mathcal{F}(\jump{\bar{u}}) \big)_{\Omega_K} - \big\langle \,1\, , \avg{\nabla u'}-\raisebox{-2pt}{\scalebox{1.2}{$\Phi$}}^{\text{\tiny I}}(\bar{u}) \big\rangle_{\hat{x}} = 0 &\qquad {\hat{x}}\in \Gamma\label{appC4}
\end{align}
Many of the formulations in \cref{tab:DGforms} do not include an additional volumetric term such that $\mathcal{F}(\cdot) = 0$. In this case, it follows from \eqref{appC4}:
\begin{align}
\quad - \big\langle \,1\, , \avg{\nabla u'}-\raisebox{-2pt}{\scalebox{1.2}{$\Phi$}}^{\text{\tiny I}}(\bar{u})\big\rangle_{\hat{x}} = 0 \quad \Rightarrow \quad \avg{\nabla u'}\Big|_{{\hat{x}}} = \raisebox{-2pt}{\scalebox{1.2}{$\Phi$}}^{\text{\tiny I}}(\bar{u})\Big|_{{\hat{x}}} &\qquad {\hat{x}}\in \Gamma\label{appC5}
\end{align}

Any coarse-scale test function $\bar{w}\in\widebar{\mathcal{V}}(0)$ may be written as the sum of discontinuous functions of the form shown in \cref{fig:appC2}, plus a continuous function. Due to linearity of $\mathcal{F}(\cdot)$, an identity similar to \cref{appC4} can be obtained for each of the discontinuous components of $\bar{w}$. Therefore, the second and third term in \cref{appC1} always cancel out. In addition, due to the nodal identity \cref{appC3}, the fourth term in \cref{appC1} vanishes as well. We thus can obtain the following element-wise volumetric identity: 
\begin{align}
-\big(\Delta \bar{w}, u' \big)_{\Omega_K} = 0 \qquad \forall\, \bar{w}\in\widebar{\mathcal{V}}(0) \label{appC2}
\end{align}
Next, we show the implications of these fine-scale models for the particular example of the symmetric interior penalty method.


\subsubsection*{The symmetric interior penalty method}
According to  \cref{tab:DGforms}, the classical symmetric interior penalty (IP) formulation of the Poisson equation reads:
\begin{align}
\begin{split}
\big(\nabla \bar{w},\nabla \bar{u}\big)_{\Omega_K} - \big\langle\jump{\bar{w}},\avg{\nabla \bar{u}}\big\rangle_{\Gamma} &- \big\langle\avg{\nabla \bar{w}} ,\jump{\bar{u}}\big\rangle_{\Gamma} + \\
&\big\langle\eta h^{-1} \jump{\bar{w}}, \jump{\bar{u}}\big\rangle_{\Gamma}  = \big(\bar{w},f\big)_{\Omega_K}  \qquad \forall \,\bar{w}\in\widebar{\mathcal{V}}(0) \label{IPWF}
\end{split}
\end{align}

We argue that \eqref{IPWF} can be retrieved from the variational coarse-scale formulation~\eqref{coarsescale3} if we choose the following implicit fine-scale interface model:
\begin{align}
\begin{cases}
\begin{alignedat}{4}
&\avg{u'}\Big|_{{\hat{x}}} &&= \raisebox{-2pt}{\scalebox{1.2}{$\Phi$}}^{\text{\tiny I}}\Big|_{{\hat{x}}} &&= \,\, 0  &\quad {\hat{x}}\in \Gamma \\ 
&\avg{\nabla u'}\Big|_{{\hat{x}}} &&= \raisebox{-1pt}{\scalebox{1.1}{$\Theta$}}^{\text{\tiny I}} (\jump{\bar{u}})\Big|_{{\hat{x}}}\, &&= -\eta h^{-1} \jump{\bar{u}}\Big|_{{\hat{x}}} &\qquad {\hat{x}}\in \Gamma
\end{alignedat}
\end{cases}\label{IPmodel}
\end{align} 
As was just shown, these expressions are pointwise identities in the one-dimensional case.
The first line in \eqref{IPmodel} states that the average of the fine-scale solution is zero across element interfaces. Therefore, it holds that $\avg{\bar{u}} = u$ at element boundaries. It is clear by intuition that this is a sound model approach, as it ensures that the coarse-scale solution does not drift away from the exact solution. 

Although the solution $\bar{u}$ centers around the exact solution on element interfaces, jumps at element interfaces are still possible. They are suppressed by the condition in the second line in \cref{IPmodel}, whose effect can be illustrated best by the following reformulation. 
We first define a distance~$d$~as:
\begin{align}
d=\frac{1}{2}h\eta^{-1} \label{dIP}
\end{align}
Then, for any interface point $ {\hat{x}}\in \Gamma$:
\begin{align}
\begin{split}
2d \, \avg{\nabla u'}\Big|_{{\hat{x}}} &= -\jump{\bar{u}}\Big|_{{\hat{x}}}= \jump{u'}\Big|_{{\hat{x}}}  \\[0.1cm]
  d \nabla u'^+ + d \nabla u'^- &= n^{+} u'^{+} + n^{-} u'^{-}   \\[0.1cm]
n^{-} u'^{+} +  d \nabla u'^+ &= n^{-} u'^{-} -  d \nabla u'^-   \\[0.1cm]
u'^{+} - d\, n^{+} \nabla u'^+ &=u'^{-} - d \, n^{-} \nabla u'^-  \\[0.1cm]
u'^{+}({\hat{x}}-d n^{+})&\approx u'^{-}({\hat{x}}-d n^{-})  
\end{split} \label{IPmodelImp}
\end{align}
In the last line, we can identify the fine-scale interface model as a first order Taylor approximation of the fine-scale solution $u'$ at the distance $d$ from the element interface.

Based in this analysis, the second line in \cref{IPmodel} states that the fine-scale solutions in neighboring elements must be approximately equal at some distance away from the interface. The two conditions \cref{IPmodel} thus control the fine-scale solution close to the interface, thereby effectively limiting the element-to-element oscillations in the coarse-scale solution. 

Accurate Taylor approximations require small distances $d$. As the free parameter $\eta$ in \cref{IPmodel} is increased, the distance $d$ approaches zero. As a consequence, the fine-scale solution is forced to zero at element interfaces and we obtain the nodally exact coarse-scale solution $\bar{u}$. This behavior is in accordance with the penalization interpretation that gave the IP method its name.
\begin{figure}[H]
\vspace{0.3cm}
\includegraphics[width=\textwidth, trim={0.5cm 0.2cm 0.0cm 0.0cm},clip]{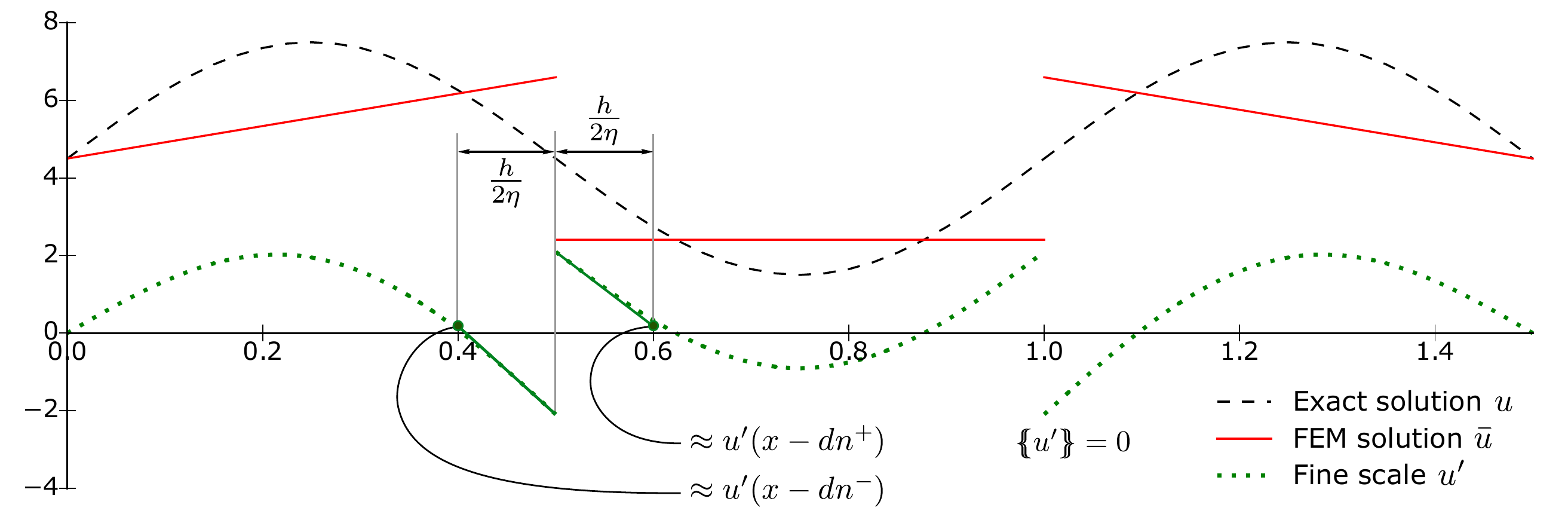}
\caption{Interior penalty method ($\eta = 2.5$): Coarse- and fine-scale solutions of a three-element linear DG discretization and graphical interpretation of the effect of the fine-scale interface model.}
\label{fig:IP}
\end{figure}
 \Cref{fig:IP} illustrates the effect of the fine-scale interface model that leads to the IP method. Using a three element linear DG discretization and a sinusoidal forcing term, we obtain the solution of the variational coarse-scale formulation \eqref{coarsescale3} with \eqref{IPmodel}. 
 We observe that the average of the discontinuous fine-scale solution $u'$ at the element interface is zero. Moreover, the first order Taylor approximations of $u'$ in the positive and negative directions from the interface at a distance $d$ are equal. 

\subsection{Higher-order discretizations of the VMS discontinuous Galerkin formulation}
\label{ssec:4.4}

Up to this point, we have restricted ourselves to DG discretizations with linear basis functions. Then the volumetric fine-scale term disappears, as it involves second derivatives of the test function. For higher-order DG discretizations, however, the volumetric term remains. 
None of the classical DG formulations shown in \cref{tab:DGforms}, however, include a volumetric term of the form $(\Delta \bar{w},\, \cdot\,)_{\Omega_K}$.
We showed in \cref{appC2} that, for formulations of the form \cref{genClassform}, the omission of the volumetric fine-scale term means it must be equal to zero. As a consequence, each of the classical DG formulations implicitly includes the following fine-scale volumetric model:
\begin{align}
\big(\Delta \bar{w}, u'\big)_{\Omega_K} = 0\label{vol0}
\end{align}
From \eqref{vol0}, it follows that for quadratic basis functions in $\bar{w}$, the element average of the fine-scale solution $u'$ is zero. For third-order basis functions, the integral of $u'$ weighted by a linear function equals zero. Generalizing this line of thought, we can write for any polynomial order $p$ higher than linear:
\begin{align}
p>1:\qquad && \int\limits_K u' \sum\limits_{n=1}^{p-1}\, c_n x^{n-1} = 0 \qquad  && \forall\, K\in \mathcal{T},\,\, \forall\, c_1, \cdots , c_{p-1} \in \mathbb{R}^{p-1}
\label{stablHO}
\end{align}
We anticipate that \eqref{stablHO} exhibits a stabilizing effect on the coarse-scale solution. The more the fine-scale solution is constrained to zero, the better the coarse-scale solution approximates the exact solution. 

\begin{figure}[H]
\includegraphics[width=\textwidth, trim={1.0cm 0.2cm 0.0cm 0.0cm},clip]{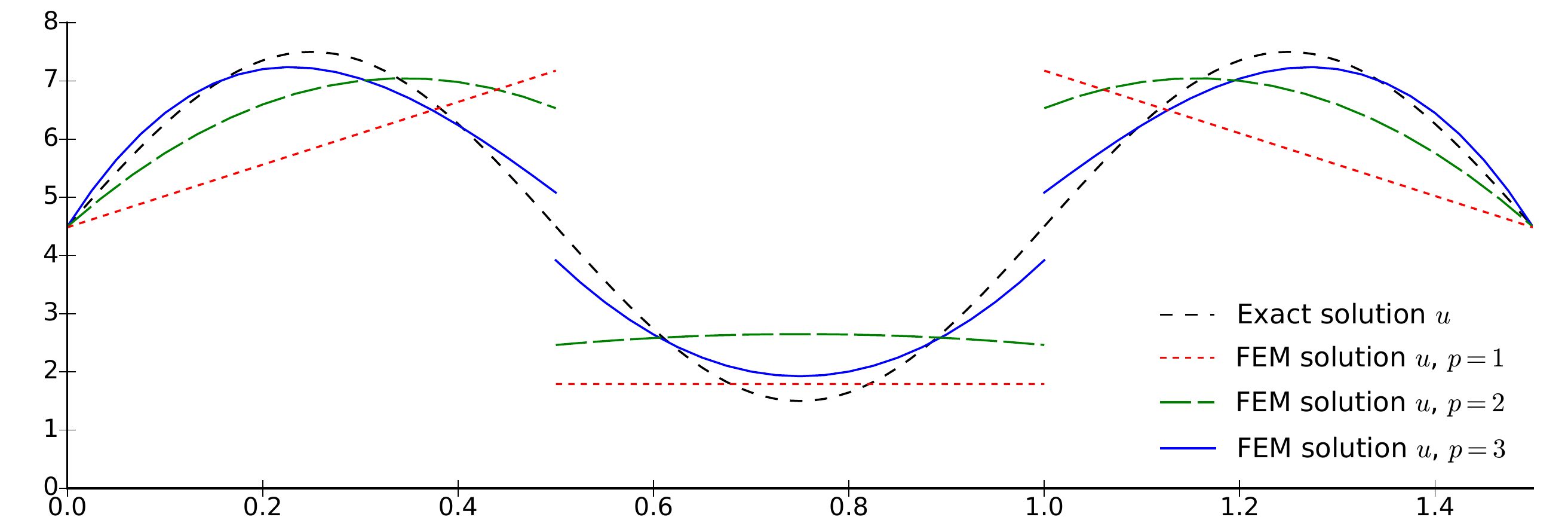}
\caption{Effect of volumetric fine-scale terms: Coarse-scale solutions of the interior penalty method with $\eta=2$. Discretized with linear, quadratic and cubic DG elements.}
\label{fig:BO2}
\end{figure}

In \Cref{fig:BO2}, we illustrate the effect of \eqref{stablHO} for higher-order discretizations of the interior penalty method and the model problem considered in the previous section. 
We verified the identities \cref{stablHO} numerically for all solutions shown in \cref{fig:BO2}. We emphasize that the solution still adheres to fine-scale interface model \eqref{IPmodel}. 

\subsection{Discretizations of the VMS discontinuous Galerkin formulation on multi-dimensional domains}
\label{ssec:4.5}
The conclusions drawn for discretizations of the VMS discontinuous Galerkin formulation on one-dimensional domains are all based on fine-scale interface models that hold pointwise at element interface nodes. On two- and three-dimensional domains, the argument made at the start of \Cref{ssec:4.3} does not hold, and pointwise identities may not be concluded from the obtained formulations. 
In this section we aim to take the first step towards developing the multiscale interpretation of DG formulations for multi-dimensional problems.

The numerical experiments in this section involve the weak enforcement of Dirichlet boundary conditions. 
We add corresponding boundary terms to the coarse-scale variational formulation \cref{coarsescale3}, leading to:
\begin{align}
 \begin{split}
&\text{Find }\bar{u}\in\widebar{\mathcal{V}}\text{ s.t.:}\\
&\quad\big(\nabla \bar{w},\nabla \bar{u}\big)_{\Omega_K}  - \big\langle\bar{w} ,\nabla \bar{u}\cdot n \big\rangle_{\partial \Omega} - \big\langle\jump{\bar{w}},\avg{\nabla \bar{u}}\big\rangle_{\Gamma}  - \big\langle\avg{\nabla \bar{w}} ,\jump{\bar{u}}\big\rangle_{\Gamma} + \\
&\qquad\qquad\quad \big(\mathcal{L}^*\!\bar{w}, u' \big)_{\Omega_K} - \big\langle\bar{w} ,\nabla u'\cdot n \big\rangle_{\partial \Omega} - \big\langle\jump{\bar{w}},\avg{\nabla u'}\big\rangle_{\Gamma}\\
&\qquad\hspace{2cm}  + \big\langle\nabla \bar{w}\cdot n, u' \big\rangle_{\partial \Omega}  + \big\langle\jump{\nabla \bar{w}}, \, \avg{u'}\big\rangle_{\Gamma}  = \big(\bar{w},f\big)_{\Omega_K}  \qquad \forall \,\bar{w}\in\widebar{\mathcal{V}} \label{coarsescale4}
\end{split}
\end{align}
Inspired by the one-dimensional case discussed earlier, we substitute the following implicit fine-scale interface model:
\begin{align}
\begin{cases}
\begin{alignedat}{4}
&\avg{u'}\Big|_{{\hat{x}}} &&= \,\, 0  \qquad&& {\hat{x}}\in \Gamma \\ 
&\avg{\nabla u'}\Big|_{{\hat{x}}} &&= -\eta h^{-1} \jump{\bar{u}}\Big|_{{\hat{x}}} \qquad&& {\hat{x}}\in \Gamma \\
&\nabla u'\cdot n \Big|_{{\hat{x}}} &&= -\eta h^{-1}\, u' \Big|_{{\hat{x}}} \qquad&& \hat{x}\in \partial \Omega
\end{alignedat}
\end{cases}\label{multiDIP}
\end{align}
Using the identity \mbox{$u'=u_D-\bar{u}$} to remove any remaining fine-scale terms on the domain boundary, we retrieve the classical IP formulation:
\begin{align}
\begin{split}
&\big(\nabla \bar{w},\nabla \bar{u}\big)_{\Omega_K} - \big\langle\jump{\bar{w}},\avg{\nabla \bar{u}}\big\rangle_{\Gamma} - \big\langle\avg{\nabla \bar{w}} ,\jump{\bar{u}}\big\rangle_{\Gamma} + \big\langle\eta h^{-1} \jump{\bar{w}}, \jump{\bar{u}}\big\rangle_{\Gamma} \\
&\hspace{1.5cm}+ \big\langle\eta h^{-1}  \bar{w} ,\bar{u} \big\rangle_{\partial \Omega} - \big\langle\nabla \bar{w}\cdot n ,\bar{u} \big\rangle_{\partial \Omega} - \big\langle\bar{w}, \nabla \bar{u}\cdot n \big\rangle_{\partial \Omega}\\
&\hspace{2.5cm} = \big(\bar{w},f\big)_{\Omega_K} -\big\langle\nabla \bar{w}\cdot n ,u_D \big\rangle_{\partial \Omega}+\big\langle\eta h^{-1}  \bar{w} ,u_D \big\rangle_{\partial \Omega} \qquad \forall \,\bar{w}\in\widebar{\mathcal{V}} \label{IPmultiD}
\end{split}
\end{align}
The formulation \cref{IPmultiD} does not satisfy the fine-scale interface model \cref{multiDIP} pointwise. Instead, it constitutes a special case of \cref{coarsescale4}, for which the fine-scale solution satisfies:
\begin{align}
\begin{alignedat}{3}
&-\big(\Delta \bar{w}, u' \big)_{\Omega_K} - \big\langle\bar{w} ,\nabla u'\cdot n \big\rangle_{\partial \Omega} - \big\langle\jump{\bar{w}},\avg{\nabla u'}\big\rangle_{\Gamma}
+ \big\langle\nabla \bar{w}\cdot n, u' \big\rangle_{\partial \Omega}  \\
&\hspace{1.7cm} + \big\langle\jump{\nabla \bar{w}}, \, \avg{u'}\big\rangle_{\Gamma}  =  \big\langle\eta h^{-1} \jump{\bar{w}}, \jump{\bar{u}}\big\rangle_{\Gamma}  - \big\langle\eta h^{-1}  \bar{w} ,u' \big\rangle_{\partial \Omega} \qquad \forall \,\bar{w}\in\widebar{\mathcal{V}} \label{IPmultiD_model}
\end{alignedat}
\end{align}
In a multi-dimensional case, this \textit{is} the fine-scale model. It forms a weakly satisfied condition that combines volumetric and interface terms.

In analogy to the one-dimensional case, the multi-dimensional VMS formulation provides tangible insights into the fine-scale behavior of classical DG formulations. For instance, by taking a test function $\bar{w}$ that equals 1 in the domain of the element $K$ and 0 outside, we find:
\begin{align}
\int\limits_{\partial K} \avg{\nabla u'}\cdot n = - \int\limits_{\partial K} \eta h^{-1} \jump{u'} \label{IP2D}
\end{align}
where the jumps and averages on the domain boundary ${\partial K}$ are defined as $\jump{u'}\equiv u'$ and $\avg{\nabla u'}\equiv \nabla u'$. 
This relation looks very much like the second line in the fine-scale interface model \eqref{IPmodel} of the one-dimensional IP method. In contrast, similar expressions for the first line of \eqref{IPmodel} and the volumetric fine-scale model of \eqref{vol0} do not hold, i.e.:
\begin{align}
\int\limits_{\partial K} \avg{u'} \neq 0 \;\;\;\;\; \text{and} \;\;\;\;\; \int\limits_{K} \,\, u' \neq 0
\label{int23}
\end{align}

\begin{figure}[ht]
\hspace{-0.6cm}\subfloat{%
\begin{tabular}{c}\vspace{-7.5cm}\\
\textbf{Model problem} \\
$\Omega = [\,0\,,\,1\,]\times[\,0\,,\,1\,]\hspace{0.75cm}$\\
$\qquad f=0 \hspace{1.72cm}\text{in } \Omega\hspace{1.28cm} $ \\
$\qquad u = 0 \hspace{1.7cm}\text{when } x_1 = 0$ \\
$\qquad u = 0 \hspace{1.7cm}\text{when } x_1 = 1$ \\
$\qquad u = \sin(\pi\,x_1) \hspace{.55cm}\text{when } x_2 = 0$ \\
$\qquad u = 0 \hspace{1.7cm}\text{when } x_2 = 1$ \\\\
\textbf{Exact solution} \\
$u(${\footnotesize $x_1,x_2$}$) =$ \\[0.1cm]
\footnotesize$\qquad \footnotesize \left(\cosh(\pi\,x_2)-\dfrac{ \cosh(\pi)}{\sinh(\pi)} \sinh(\pi\,x_2) \right)\,\sin(\pi\,x_1)$    \\ \\
\textbf{Numerical implementation}\\
Mesh: $18$ triangular elements  \\
IP method with $\eta_{_\Gamma} = 3$ and $\eta_{_{\partial \Omega}} = 8$
\end{tabular}
}\hspace{-0.2cm}
\subfloat{%
  \includegraphics[width=0.5\textwidth,trim={10cm 5.cm 7.5cm 1cm},clip]{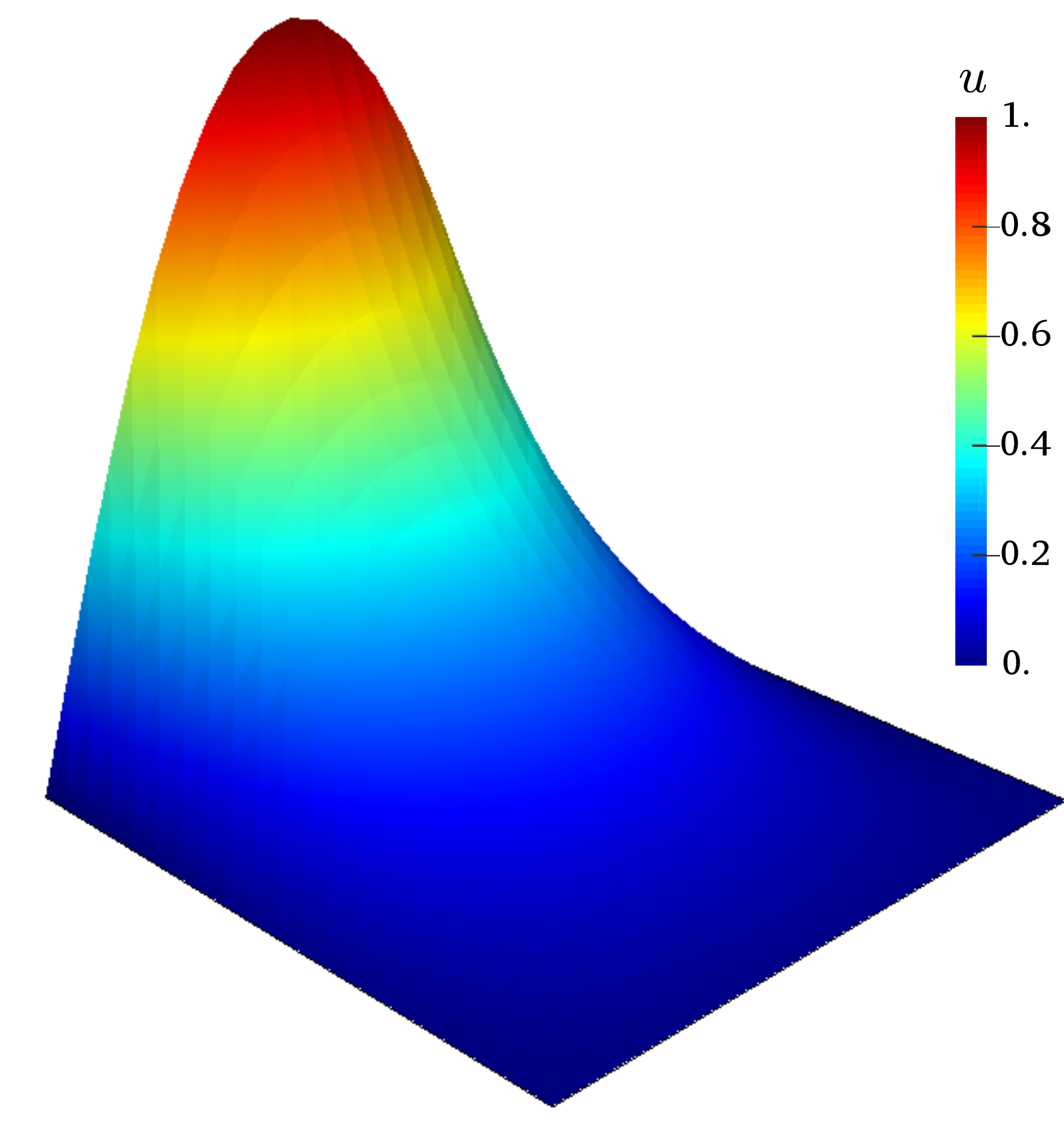}
}
  \caption{Exact solution of the 2-dimensional Laplace model problem.}%
  \label{fig:exp2D}
\end{figure}

\begin{figure}[H]
\centering
\subfloat[Coarse-scale solution $\bar{u}$.] {\includegraphics[width=0.48\textwidth,trim={10cm 5.cm 7.5cm 1cm},clip]{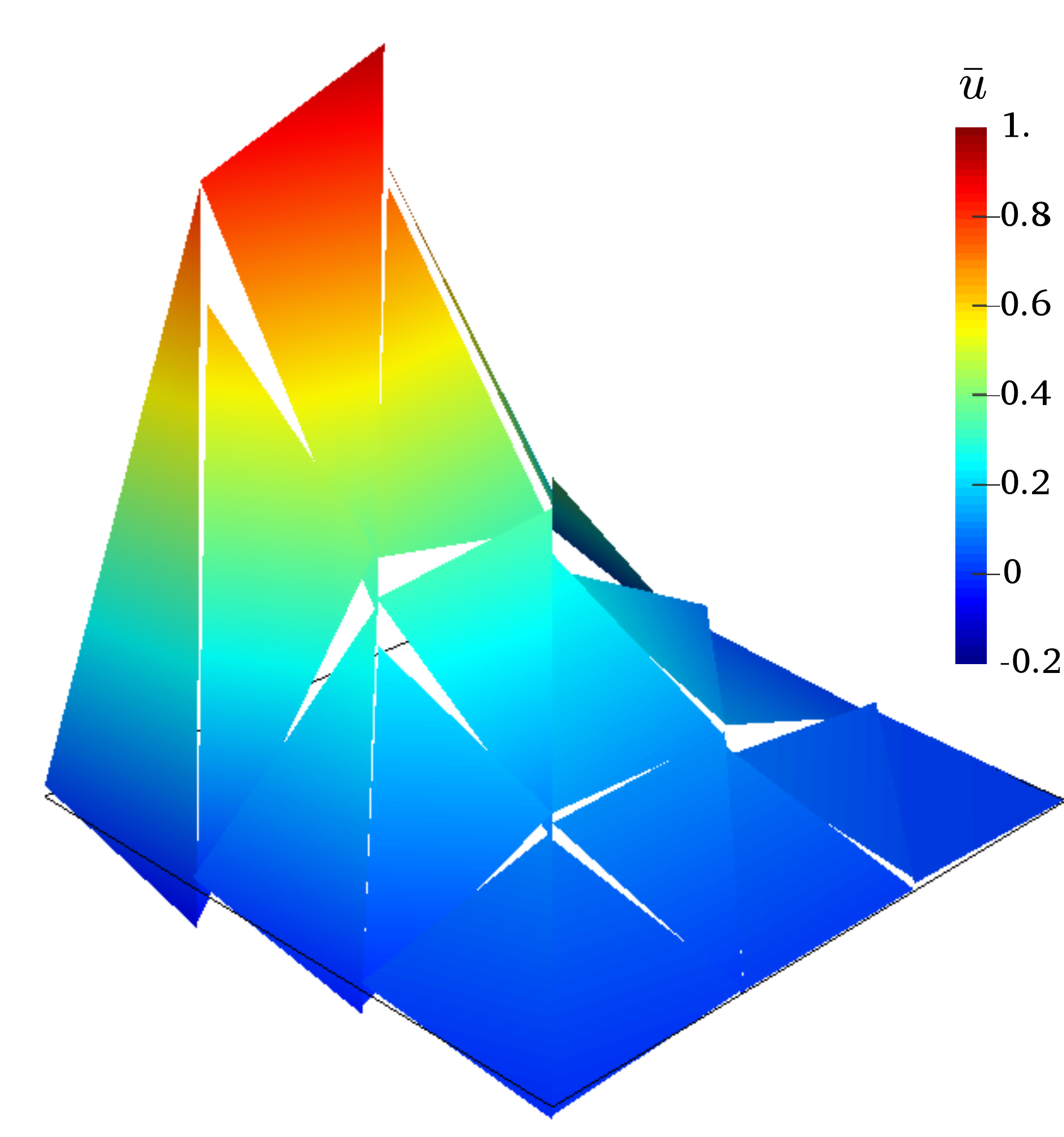} }  \hspace{0.1cm}
\subfloat[Fine-scale solution $u'$.] {\includegraphics[width=0.48\textwidth,trim={10cm 5.cm 7.5cm 1cm},clip]{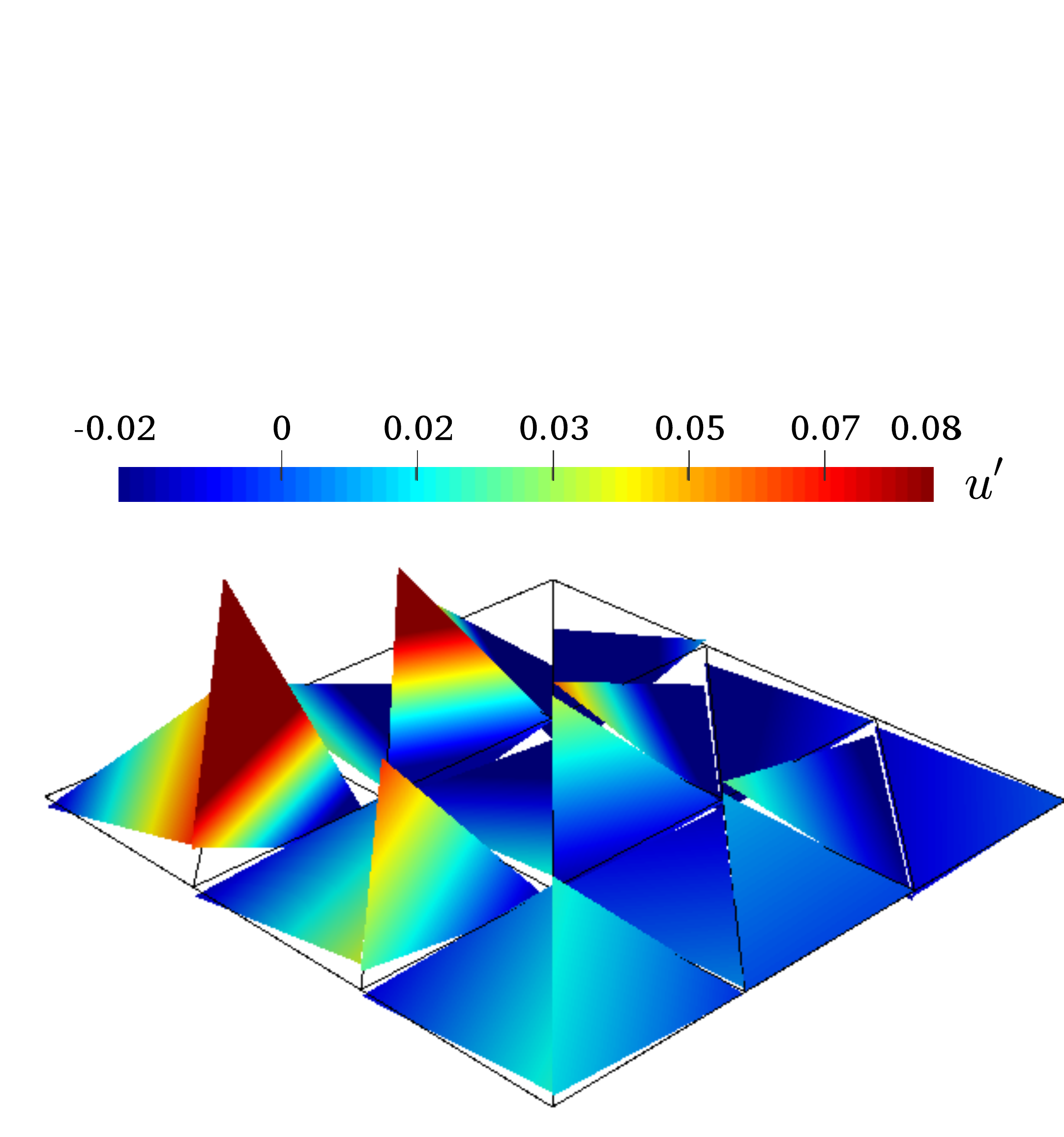} }
\caption{Interior penalty formulation, discretized with 18 linear DG elements ($p=1$).}
\label{fig:2D1p}
\end{figure}

\begin{figure}[H]
\vspace{1.2cm}
\centering
\subfloat[Coarse-scale solution $\bar{u}$.] {\includegraphics[width=0.48\textwidth,trim={10cm 5.cm 7.5cm 1cm},clip]{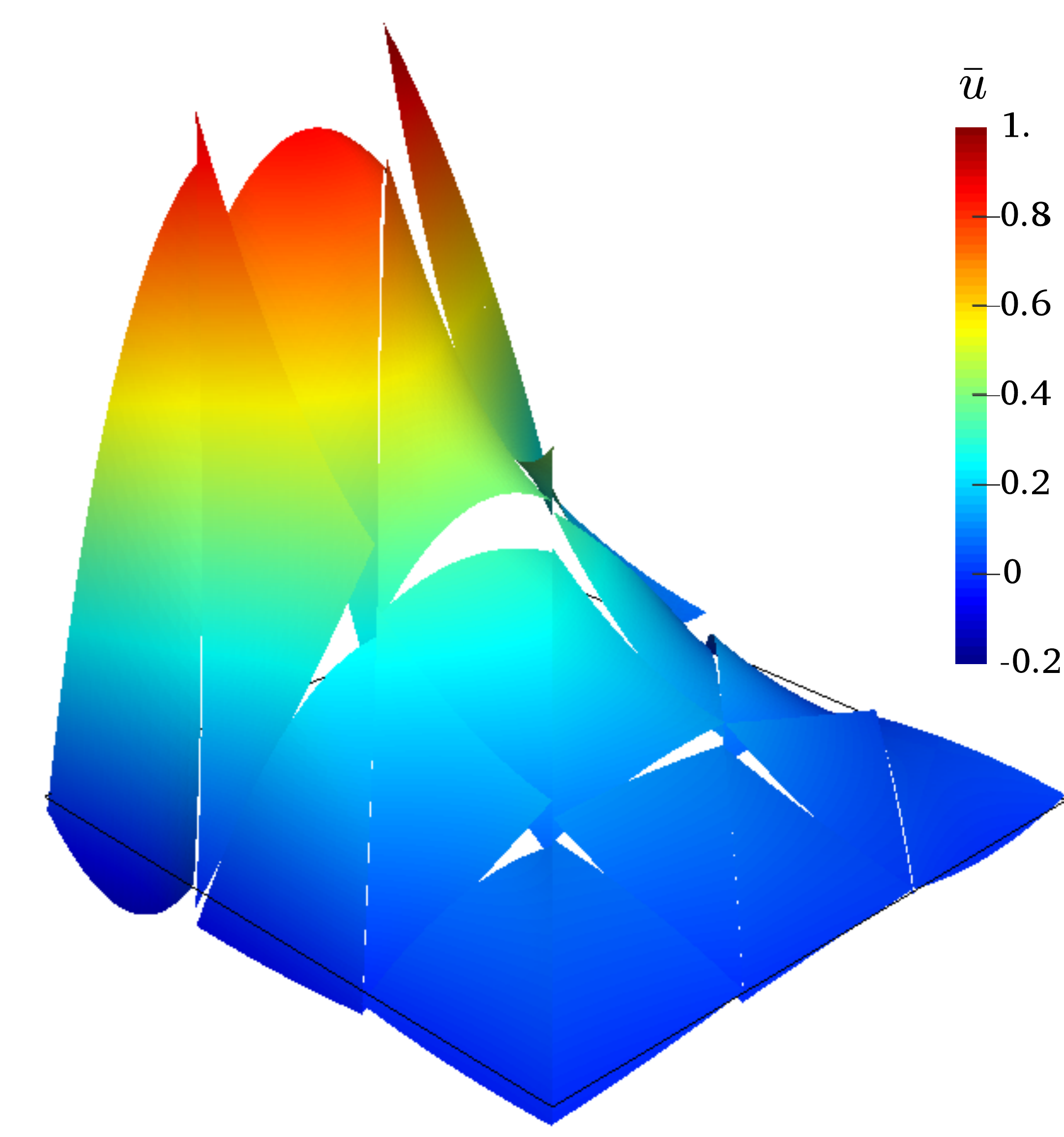} }  \hspace{0.1cm}
\subfloat[Fine-scale solution $u'$.] {\includegraphics[width=0.48\textwidth,trim={10cm 5.cm 7.5cm 1cm},clip]{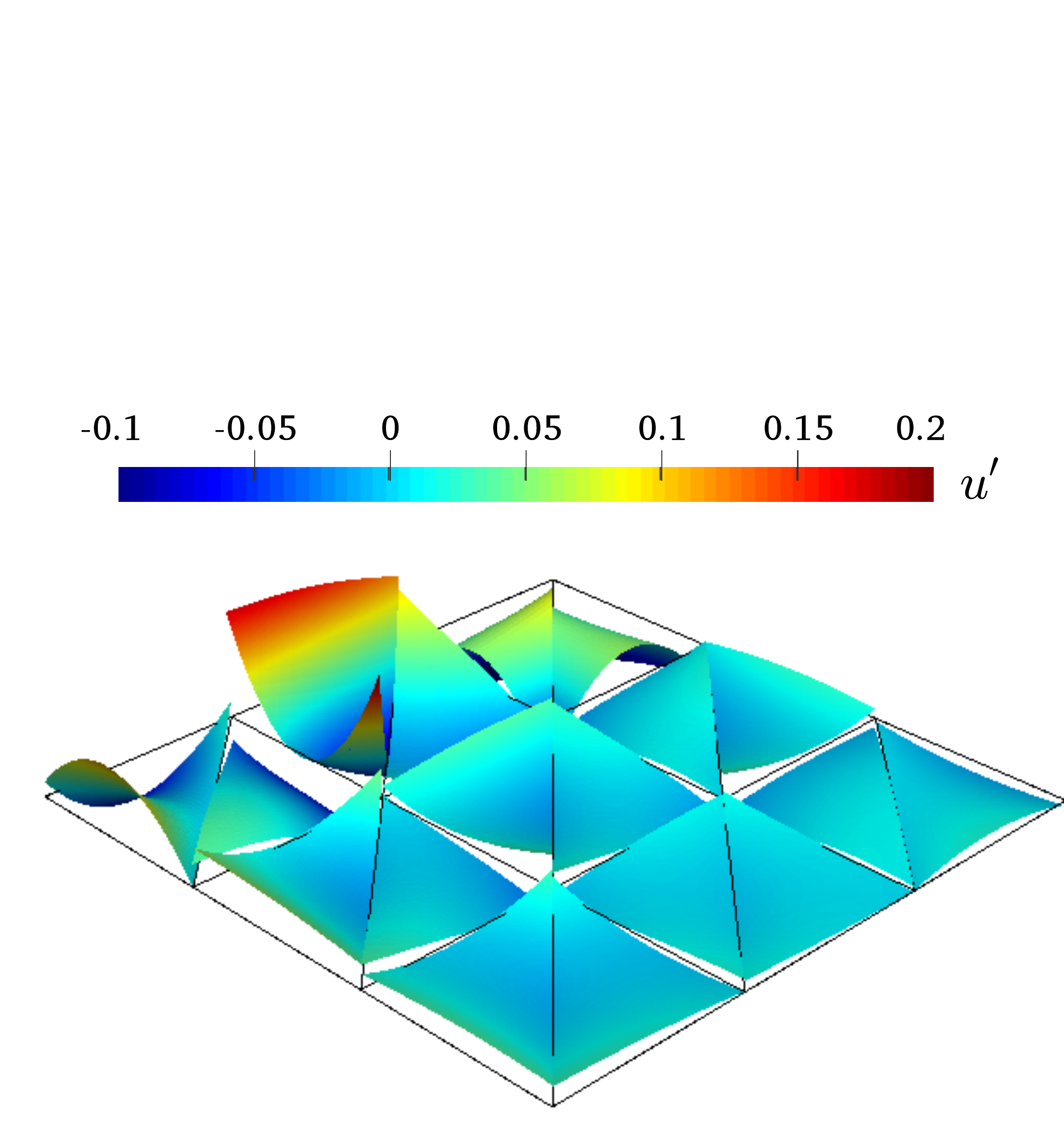} }
\caption{Interior penalty formulation, discretized with 18 quadratic DG elements ($p=2$).}
\label{fig:2D2p}
\end{figure}

\begin{figure}[H]
\centering
\subfloat[Coarse-scale solution $\bar{u}$.] {\includegraphics[width=0.48\textwidth,trim={10cm 5.cm 7.5cm 1cm},clip]{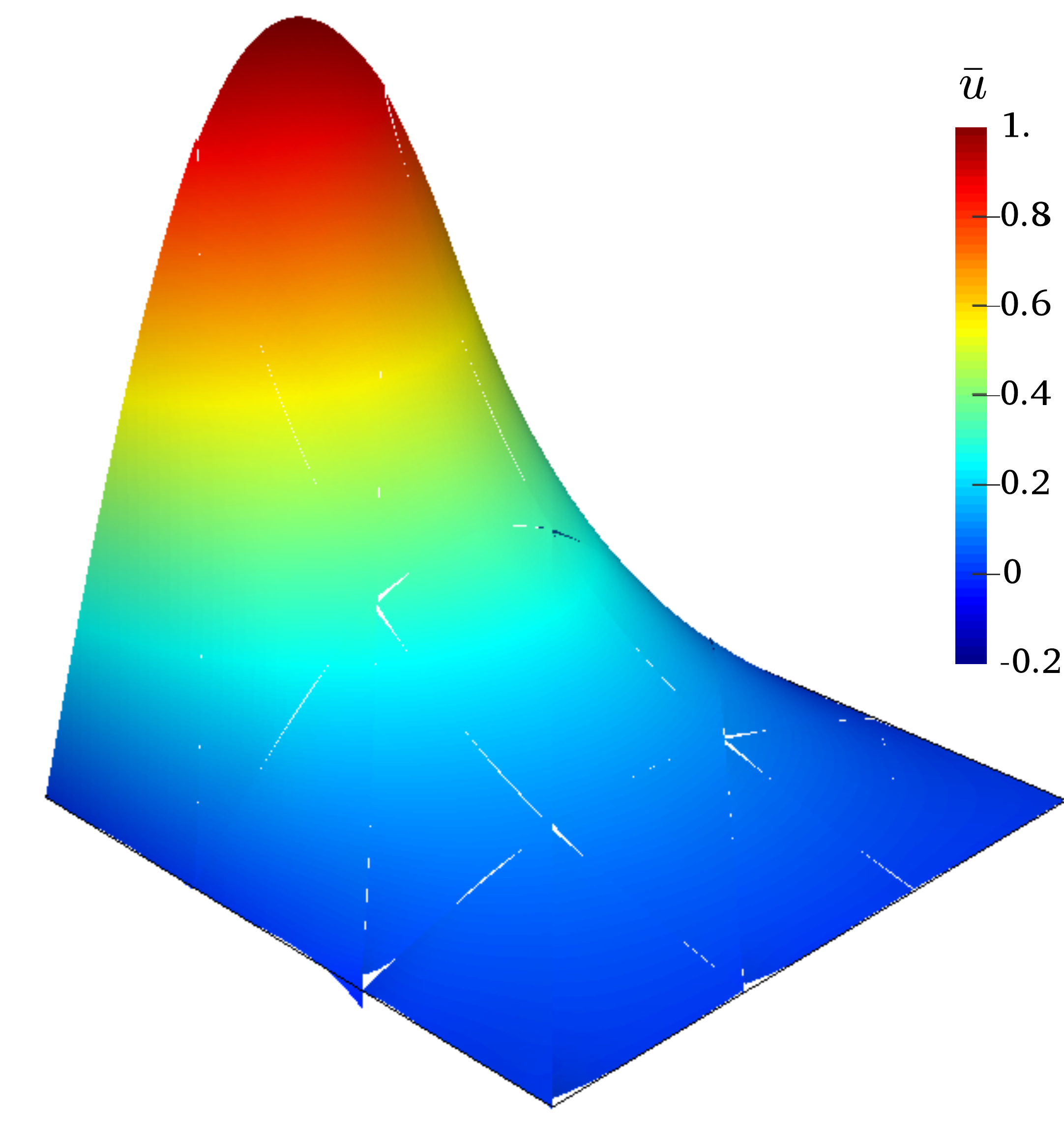} }  \hspace{0.1cm}
\subfloat[Fine-scale solution $u'$.] {\includegraphics[width=0.48\textwidth,trim={10cm 5.cm 7.5cm 1cm},clip]{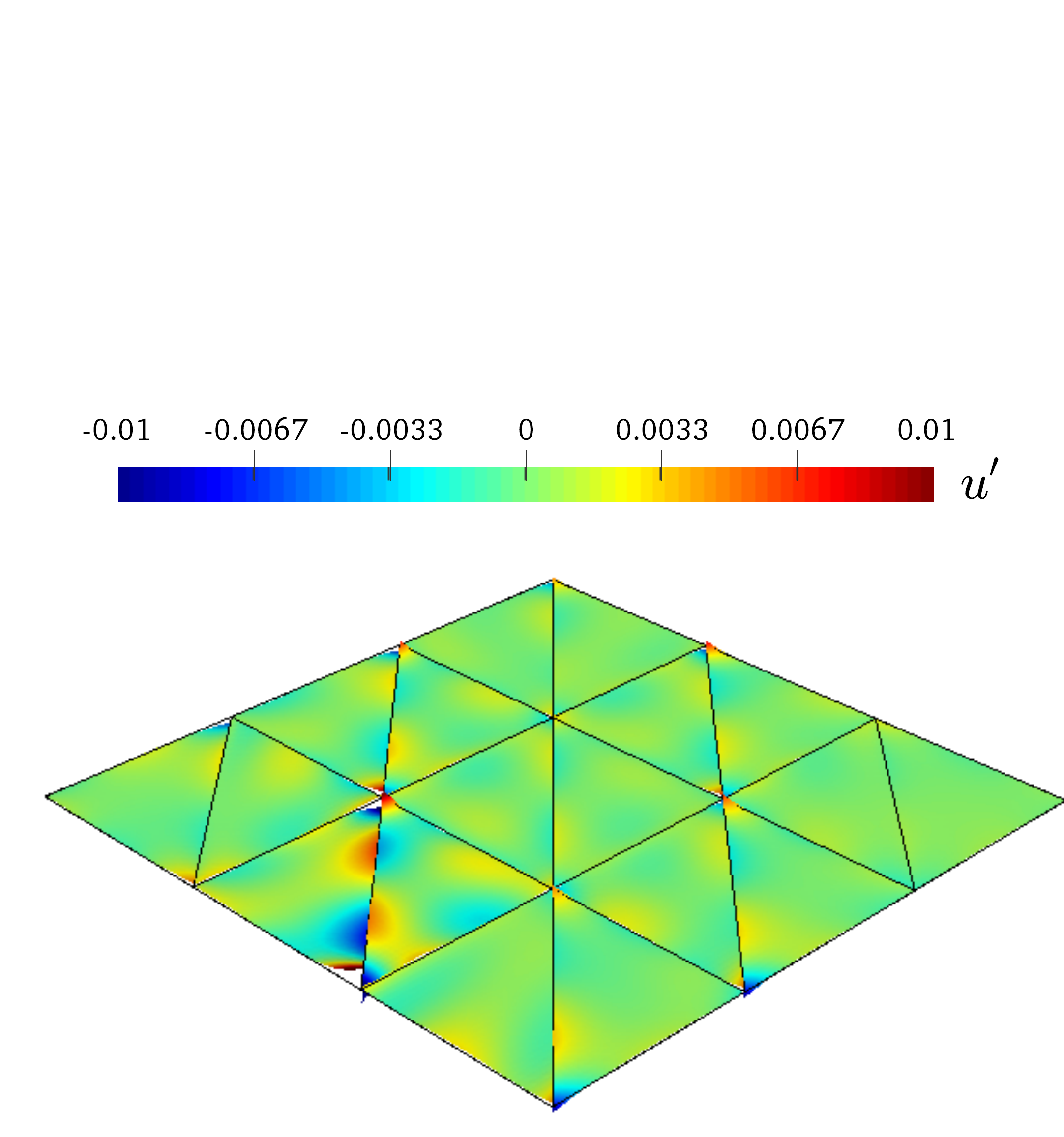} }
\caption{Interior penalty formulation, discretized with 18 quartic DG elements ($p=4$).}
\label{fig:2D4p}
\end{figure}

We investigate the behavior of the relations \eqref{int23} with a numerical experiment. To this end, we consider the 2-D Laplace problem defined in \cref{fig:exp2D}. We use a DG discretization that consists of 18 triangular elements. We increase the polynomial order of the basis functions for a number of consecutive experiments. 
\Cref{fig:2D1p,fig:2D2p,fig:2D4p} plot the coarse- and fine-scale solutions obtained with linear, quadratic and quartic basis functions, respectively. We list the results of the integral expression \cref{IP2D,int23} in \Cref{tab:NumRes2D}.

\vspace{0cm}

\vspace{0.5cm}
\newpage

\begin{table}[H]
\centering
\caption{Integrals \cref{IP2D,int23} evaluated with solutions obtained with the IP method for the problem defined in \cref{fig:exp2D}. We show results for an 18-element DG discretization with basis functions of polynomial order $p=1$ to $p=6$.}
\label{tab:NumRes2D}
\begin{tabular}{lccc}
\hline\hline \\[-0.3cm]
      &   $\,\, \displaystyle\max\limits_{K\in\mathcal T}\Big|\! \int\limits_{\partial K}\!\! \avg{\nabla u'}\!+\!\eta h^{-1}\jump{\bar{u}} \,\Big|\,\, $ & $\,\, \raisebox{-1.2pt}{$\frac{1}{\raisebox{-3pt}{\small\#Els}}$}\displaystyle\sum\limits_{K\in\mathcal T}\Big| \! \int\limits_{\partial K} \avg{u'} \,\Big|\,\,$ & $\,\,\raisebox{-1.2pt}{$\frac{1}{\raisebox{-3pt}{\small\#Els}}$}\displaystyle\sum\limits_{K\in\mathcal T}\Big|  \int\limits_{K} u' \,\,\Big|\,\,$ \\[0.6cm] \hline\\[-0.15cm]
$p=1$ & 0                     & $ 2.39 \cdot 10^{-3}$                    & $2.93 \cdot 10^{-4}$                                             \\[0.2cm]
$p=2$ & 0                     & $ 1.81 \cdot 10^{-3}$                    & $1.98 \cdot 10^{-4}$                                             \\[0.2cm]
$p=3$ & 0                     & $ 6.83 \cdot 10^{-5}$                    & $2.78 \cdot 10^{-6}$                                             \\[0.2cm]
$p=4$ & 0                     & $ 2.22 \cdot 10^{-5}$                    & $2.04 \cdot 10^{-6}$                                             \\[0.2cm]
$p=5$ & 0                     & $ 7.89 \cdot 10^{-7}$                    & $7.53 \cdot 10^{-9}$                                             \\[0.2cm]
$p=6$ & 0                     & $ 3.31 \cdot 10^{-7}$                    & $3.59 \cdot 10^{-9}$                                             \\[0.2cm]
\hline \hline                              
\end{tabular}
\end{table}

The results in the first column of \cref{tab:NumRes2D} verify that the fine-scale interface model \cref{IP2D} is exactly satisfied. We observe in the second and third columns that the integrals \cref{int23} are not zero, but they approach zero with increasing polynomial order of the basis functions. Note that there is an even-odd phenomenon with respect to the order of the basis functions. The results of the integrals reduce by several orders of magnitude when the polynomial order is increased from an even to an odd value, but stay practically constant when the order is increased from odd to even. 
This interesting convergence behavior 
warrants a more thorough investigation of multi-dimensional fine-scale models in the future. 

\section{A residual-based VMS discontinuous Galerkin method for steady advection-diffusion problems}
\label{sec:5}

The discontinuous Galerkin method naturally incorporates stable upwind numerical flux formulations. It has therefore established itself as an effective tool for discretizing boundary value problems based on advection-type PDEs. In this section, 
we show that the VMS framework derived in the previous sections also holds for advection-type problems. In particular, we explore an interpretation of upwind numerical fluxes in a VMS context.
To this end, we consider 
the following problem:
\begin{align}
\begin{cases}
\begin{alignedat}{4}
a\cdot \nabla u - \nu \Delta u &= f &&\qquad \text{in } \Omega \subset \mathbb{R}^d \\
u &= u_D &&\qquad \text{on } \partial \Omega \label{AdDifProb}
\end{alignedat}
\end{cases}
\end{align}
The parameters $a$ and $\nu$ are assumed constant in our analysis.

\subsection{Variational multiscale formulation}
\label{ssec:5.1}
Following the procedures described in \Cref{sec:2}, we can derive the variational coarse-scale formulation as:
\begin{align}
\begin{alignedat}{2}
&\text{Find } \bar{u}\in \widebar{\mathcal{V}}(u_D) \text{ s.t.:}\\
&\quad-\big(a\cdot \nabla \bar{w},\bar{u}\big)_{\Omega_K}\!+\big(\nabla \bar{w},\nu \nabla \bar{u}\big)_{\Omega_K}\!+\sum\limits_{K\in \mathcal{T}}\Big[ \big\langle\bar{w},a\cdot n \,\bar{u}\big\rangle_{\partial K} \!- \big\langle\bar{w},\nu \, \nabla \bar{u}\cdot n \big\rangle_{\partial K} \Big] + \\[-0.1cm]
&\quad \hspace{2.3cm} \big(\mathcal{L}^*\!\bar{w}, u' \big)_{\Omega_K}\! + k(\bar{w},u';\Gamma) = \big(\bar{w},f\big)_{\Omega_K}  \qquad \forall \,\bar{w}\in\widebar{\mathcal{V}}(0) \label{coarsescale1AD} 
\end{alignedat}
\end{align}
where the following holds:
\begin{align}
&\mathcal{L}^* = -a\cdot \nabla - \nu\Delta \\[0.1cm]
&k(\bar{w},u';\Gamma)  = \sum\limits_{K\in \mathcal{T}}\Big[  \big\langle \bar{w},a\cdot n \,u'\big\rangle_{\partial K}  -\big\langle \bar{w},\nu\,\nabla u'\cdot n\big\rangle_{\partial K} + \big\langle\nabla \bar{w}\cdot n ,\nu\, u'\big\rangle_{\partial K} \Big] \label{kaddif}
\end{align}

We find the globally coupled formulation by manipulation of the fine-scale element boundary terms. We use \cref{step4Key} in \cref{kaddif} and substitute the result into \cref{coarsescale1AD} to obtain the following variational coarse-scale formulation:
\begin{align}
\begin{alignedat}{2}
&\text{Find } \bar{u}\in \widebar{\mathcal{V}}(u_D) \text{ s.t.:}\\
&\quad -\big(a\cdot \nabla \bar{w},\bar{u} \big)_{\Omega_K} + \big\langle\jump{\bar{w}}\cdot a,\avg{ \bar{u}} \big\rangle_{\Gamma} + \big\langle\jump{\bar{w}}\cdot a,\avg{u'} \big\rangle_{\Gamma} +\\ 
&\quad \hspace{0.4cm} \big(\nabla \bar{w},\nu\nabla \bar{u} \big)_{\Omega_K} - \big\langle\jump{\bar{w}},\nu \, \avg{\nabla \bar{u}}\big\rangle_{\Gamma} - \big\langle\avg{\nabla \bar{w}} ,\nu \, \jump{\bar{u}}\big\rangle_{\Gamma} + \\
&\quad \big(\mathcal{L}^*\!\bar{w}, u' \big)_{\Omega_K}    - \big\langle\jump{\bar{w}},\nu \, \avg{\nabla u'} \big\rangle_{\Gamma} + \big\langle\jump{\nabla \bar{w}}, \,\nu \,  \avg{u'} \big\rangle_{\Gamma}  = \big(\bar{w},f \big)_{\Omega_K}  \quad \forall \,\bar{w}\in\widebar{\mathcal{V}}(0)\label{coarsescale2AD} 
\end{alignedat}
\end{align}
The collection of the terms associated with the diffusion operator has been discussed in \Cref{ssec:4.1}. Here, we apply a similar approach that results in the two advection terms in \eqref{coarsescale2AD}. In particular, we obtain the term $\langle\jump{\bar{w}}\cdot a,\avg{ \bar{u}}\rangle_{\Gamma}$ by a manipulation similar to \cref{collection}.
We note that in \eqref{coarsescale2AD}
the volumetric fine-scale term does not disappear, 
because the adjoint differential operator $\mathcal{L}^*$ in \eqref{kaddif} includes a first derivative term. Therefore, $\mathcal{L}^*\!\bar{w}$ has non-zero values, also for linear basis functions. 

The fine-scale solution $u'$ in the element interior is obtained as described in \Cref{sec:3}. Assuming constants for $f$, $a$ and $\nu$ in the one-dimensional problem \cref{AdDifProb}, we find: 
\begin{align}
\begin{split}
\big(\mathcal{L}^*\!\bar{w}, u' \big)_{\Omega_K} =  & \sum\limits_{k\in\mathcal{T}}\Big[ \big(-a\cdot \nabla \bar{w}\,,\, \tau\, (f-a\cdot\nabla \bar{u})\, \big)_{K} + \\[-0.3cm]
&\hspace{1.2cm}\big(-a\cdot \nabla \bar{w}, \gamma_0 \nu \,n_{j}\, u'_j \big)_{K} + \big(-a\cdot \nabla \bar{w}, \gamma_1 \nu \,n_{j+1}\,u'_{j+1} \big)_{K} \Big] \label{taugammagamma}
\end{split}
\end{align}
with the average Green's function quantities:
\begin{align}
\tau\, = \frac{h}{2\,a}-\frac{\nu}{a^2}+\frac{h}{a \left( e^{\frac{a}{\nu}h}-1 \right)}; \qquad
\gamma_0 =  \frac{\nu - a\, h - \nu e^{-\frac{a}{\nu} h}}{a\, h\, \nu \left(e^{-\frac{a}{\nu} h}-1\right)};  \qquad
\gamma_1 = \frac{\nu + a\, h - \nu e^{\frac{a}{\nu} h}}{a\, h\, \nu \left(e^{\frac{a}{\nu} h}-1\right)} \label{tgg}
\end{align}
For a detailed derivation of the expressions in \eqref{tgg}, the interested reader is referred to \Cref{App:A}. 
The vectors $n_j$ and $n_{j+1}$ in \cref{taugammagamma} are the left and right normals to each of the two boundaries of a 1-D element. The fine-scale boundary values $u'_j$ and $u'_{j+1}$ correspond to the left and right boundary of each 1-D element. 

To close the formulation, we require expressions for the fine-scale element boundary values $u'_j$ and $u'_{j+1}$ in the volumetric fine-scale model. They can again be based on explicit or implicit models. 
To find an implicit model, they have to be related to the coarse-scale solution $\bar{u}$ in some way. In \Cref{ssec:4.3}, we found that  on a 1-dimensional domain, many existing formulations implicitly enforce:
\begin{align}
\quad \avg{u'}\Big|_{\hat{x}}  = \raisebox{-2pt}{\scalebox{1.2}{$\Phi$}}^{\text{\tiny{I}}}\Big|_{\hat{x}} = 0 \hspace{4.cm}\hat{x}\in\Gamma
\end{align}
Substitution thereof into \cref{step4Key} yields the fine-scale boundary values:
\begin{align}
\quad u'^{\pm}\Big|_{\hat{x}} = \frac{1}{2}\jump{u'}\Big|_{\hat{x}}\cdot n^\pm = -\frac{1}{2}\jump{\bar{u}}\Big|_{\hat{x}}\cdot n^\pm  \hspace{1.3cm}\hat{x}\in\Gamma
\end{align}
The implicit volumetric fine-scale model finally follows as:
\begin{align}
\begin{split}
\big(\mathcal{L}^*\!\bar{w}, u' \big)_{\Omega_K} = &  \sum\limits_{K\in\mathcal{T}}\Big[  \big(-a\cdot \nabla \bar{w} \,,\, \tau\,(f-a\cdot\nabla \bar{u})\, \big)_{K} \,+ \\[-0.3cm]
&\hspace{1.2cm}\big(a\cdot \nabla \bar{w}\, , \, \gamma_0 \nu \frac{1}{2}\jump{\bar{u}}\Big|_{x_{j}} \big)_{K} + \big(a\cdot \nabla \bar{w}\, , \, \gamma_1 \nu \frac{1}{2}\jump{\bar{u}}\Big|_{x_{j+1}} \big)_{K} \Big] \label{taugammagamma2}
\end{split}
\end{align}

\subsection{Numerical experiments with linear basis functions in 1-D}
\label{ssec:5.2}

We perform two numerical experiments to investigate the effect of the volumetric fine-scale term. We recall that the volumetric fine-scale term did not play a role in our previous numerical experiments for the Poisson problem in \Cref{ssec:4.2}, where it canceled for DG discretizations with linear basis functions. 

\begin{figure}[ht]
\hspace{-0.2cm}
\subfloat{%
\begin{tabular}{c}\vspace{-6cm}\\
\textbf{Model problem} \\
$f=6$, $a=0.5$, $\nu=0.15$ \\
$x_0 = 0$, $x_1 = 1$\\
$u_0 = 0$, $u_1 = 2$ \\\\
\textbf{Exact solution} \\
$u(x) = -\frac{10}{e^{\frac{1}{0.3}}-1}\left(e^{\frac{x}{0.3}}-1\right)+12x$    \\ \\
\textbf{Fine-scale model}\\
For ${\hat{x}}\in \Gamma$:  \\
$\avg{u'}\Big|_{{\hat{x}}}= 0$\\
$\avg{\nabla u'}\Big|_{{\hat{x}}}\!\!=\! \frac{\text{d} u}{\text{d}x}\Big|_{{\hat{x}}}\!\!-\!\frac{u({\hat{x}}+h)-u({\hat{x}}-h)}{2h}$ \\
$u'({\hat{x}}) = 0$ \\     
\end{tabular}
}\hspace{-0.3cm}
\subfloat{%
  \includegraphics[width=0.6\textwidth]{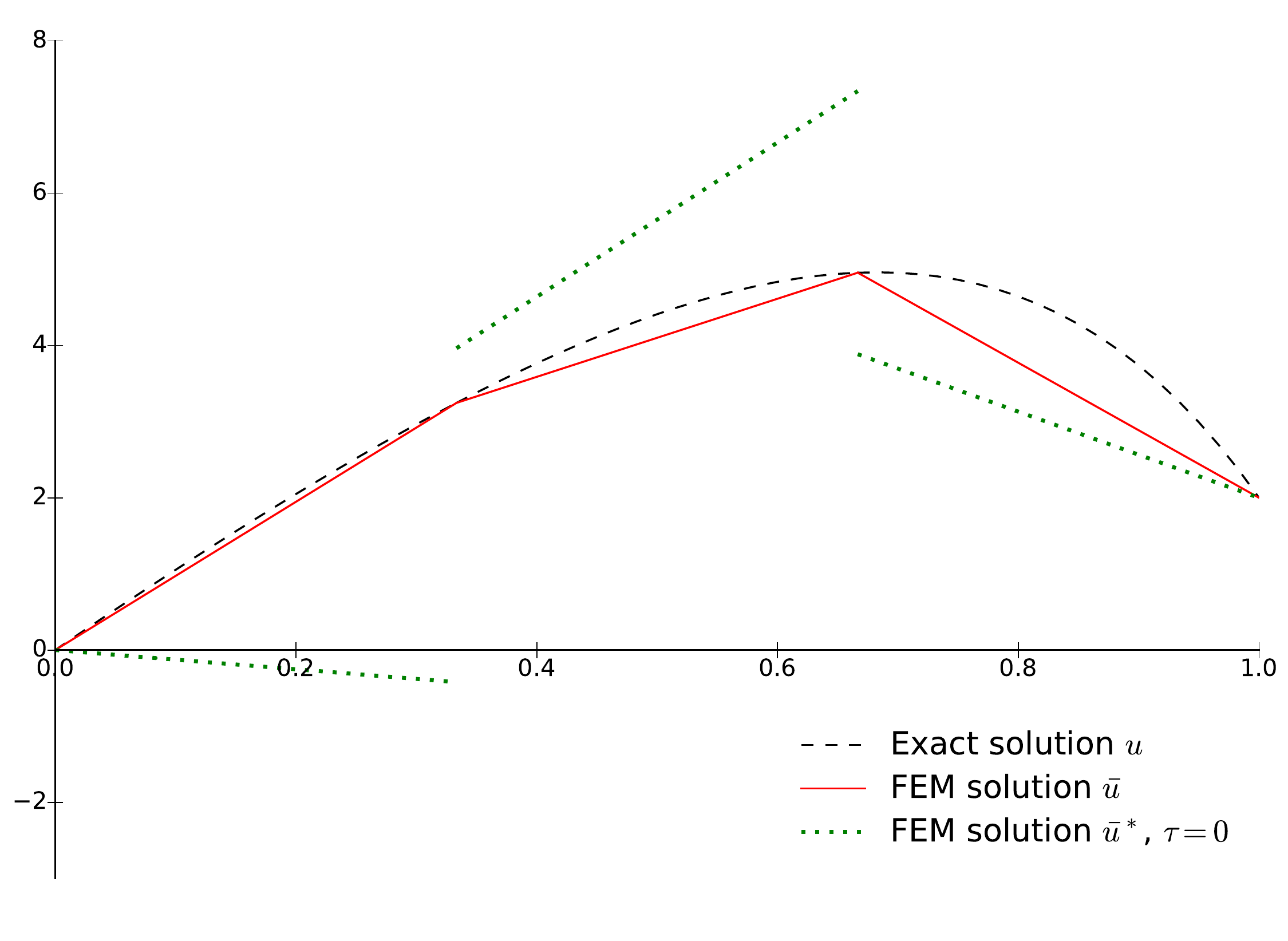}%
}
  \caption{Numerical experiment 1: the nodally exact coarse-scale solution $\bar{u}$ is obtained with an explicit formulation for the fine-scale terms.}%
  \label{fig:AdDifexp1}
\end{figure}

In the first numerical experiment, we consider a nodally exact $H^1$ projection of the exact solution for the problem defined in \Cref{fig:AdDifexp1}. In this case, all nodal values for $u'$ and $\avg{u'}$ are zero. This includes those corresponding to $\gamma_0$ and $\gamma_1$ in \cref{taugammagamma}. The fine-scale volumetric term therefore only contains the contribution associated with~$\tau$. In this case, we treat the element interface term $\avg{\nabla u'}$ explicitly, in the same fashion described for the fine-scale interface model of the Poisson problem in \Cref{ssec:4.2}.

\vspace{0.3cm}

Using a thee-element DG discretization with linear basis functions, we find the coarse-scale solution shown in \Cref{fig:AdDifexp1}. The results verify that the $H^1$ projection is retrieved. In addition, they show that the volumetric fine-scale term is essential. We observe that when the volumetric fine-scale term is canceled by setting $\tau=0$, the coarse-scale solution is far away from the exact solution. Small changes in the problem parameter $\nu$ lead to largely varying solutions, which is a well-known behavior of finite element approximations in advection-diffusion problems. 

In the second numerical experiment, we consider 
an $L^2$ projection of the exact solution, which no longer yields a nodally exact coarse-scale solution $\bar{u}$. Therefore, the fine-scale boundary values of the volumetric fine-scale model cannot be omitted here. To obtain explicit fine-scale values, we follow the same strategy as in \Cref{ssec:4.2} and precompute the correct $L^2$ projection.

\begin{figure}[t!]
\hspace{-0.2cm}
\subfloat{%
\begin{tabular}{c}\vspace{-6cm}\\
\textbf{Model problem} \\
$f=6$, $a=0.5$, $\nu=0.15$ \\
$x_0 = 0$, $x_1 = 1$\\
$u_0 = 0$, $u_1 = 2$ \\\\
\textbf{Exact solution} \\
$u(x) \!=\! -\frac{10}{e^{\frac{1}{0.3}} - 1}\left(e^{\frac{x}{0.3}}\!-\!1\right)\!+\!12x$    \\ \\
\textbf{Fine-scale model}\\
For ${\hat{x}}\in \Gamma$:  \\
\small$\avg{u'}\Big|_{{\hat{x}}}\!\!=\!u({\hat{x}}) \!-\! \avg{\mathcal{P}_{\! L^2}u}\Big|_{{\hat{x}}}$\\
\small$\avg{\nabla u'}\Big|_{{\hat{x}}}\!\!=\!\nabla u({\hat{x}}) \!-\! \avg{\nabla \,\mathcal{P}_{\! L^2}u }\Big|_{{\hat{x}}}$        \\
\small$u'\big(\limeps\,\, {\hat{x}}\pm\epsilon\big)\! \!=\! u({\hat{x}}) \!-\! \mathcal{P}_{\! L^2}u \big(\limeps\,\, {\hat{x}}\pm\epsilon\big)\Big.$\\    
\end{tabular}
}\hspace{-0.4cm}
\subfloat{%
  \includegraphics[width=0.59\textwidth,trim={0.4cm 0.cm 0.0cm 0cm},clip]{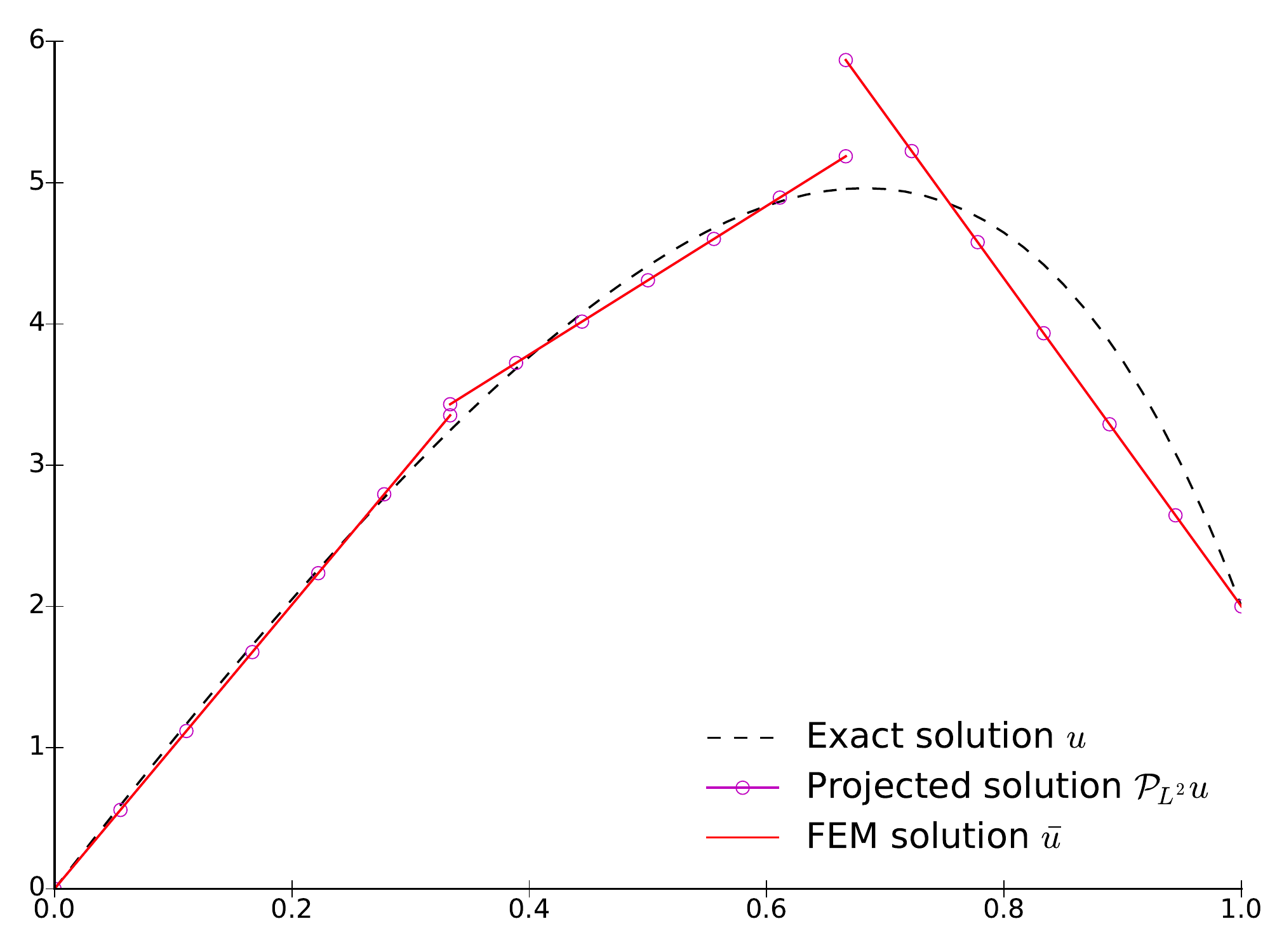}%
}
  \caption{Numerical experiment 2: the $L^2$ projection of the exact solution $u$ is obtained with an explicit formulation for the fine-scale terms.}%
  \label{fig:AdDifexp3}
\end{figure}

\Cref{fig:AdDifexp3} illustrates the result obtained with a three-element linear DG discretization. We observe that the $L^2$ projection of the exact solution is indeed retrieved by explicit substitution of the fine-scale terms. This indicates the validity of the fine-scale terms in the variational coarse-scale formulation \cref{coarsescale2AD} and the average Green's function identities \cref{taugammagamma}.

\subsection{A VMS interpretation of upwind flux formulations}
\label{ssec:5.3}

Upwind flux formulations have been the backbone of discontinuous Galerkin methods since their introduction in 1973~\cite{Reed1973}. In this subsection, we investigate upwind flux formulations from a multiscale perspective. In addition, we derive the fine-scale interface model that corresponds to the classical upwind strategy. 


The physical interpretation of advection states that information is propagated in the direction of the velocity vector $a$. In an upwinding method, the element boundary data of some element in the mesh relates to the upstream element. Inspired by this notion, we assume that the flux term obtained from integration by parts of the coarse-scale solution should have the form:
\begin{align}
\big\langle\jump{\bar{w}}\cdot a,\bar{u}(\limeps\,\, x-\epsilon a)\big\rangle_{\Gamma} \label{fluxterm}
\end{align}
where we obtain jumps in $\bar{w}$, because $\bar{u}(\limeps\,\, x-\epsilon a)$ is single-valued on the element interface. 
We obtain a typical upwind formulation by substituting \cref{exacttoupwind} into \cref{coarsescale2AD}. 
\begin{align}
\quad \big\langle\jump{\bar{w}}\cdot a,\avg{\bar{u}}\big\rangle_{\Gamma} + \big\langle\jump{\bar{w}}\cdot a,\avg{u'} \big\rangle_{\Gamma}  = \big\langle\jump{\bar{w}}\cdot a,\bar{u}(\limeps\,\, x-\epsilon a)\big\rangle_{\Gamma} \label{exacttoupwind}
\end{align}

\begin{figure}[t!]
\hspace{-0.2cm}\subfloat{%
\begin{tabular}{c}\vspace{-6cm}\\
\textbf{Model problem} \\
$f=6$, $a=0.5$, $\nu=0.001$ \\
$x_0 = 0$, $x_1 = 1$\\
$u_0 = 0$, $u_1 = 2$ \\\\
\textbf{Exact solution} \\
$u(x) = -\frac{10}{e^{500}-1}\left(e^{500x}-1\right)+12x$    \\ \\
\textbf{Fine-scale model}\\
For ${\hat{x}}\in \Gamma$:  \\
$\avg{u'}\Big|_{{\hat{x}}}= 0$\\
$\avg{\nabla u'}\Big|_{{\hat{x}}}\!\!=\! \frac{\text{d} u}{\text{d}x}\Big|_{{\hat{x}}}\!\!-\!\frac{u({\hat{x}}+h)-u({\hat{x}}-h)}{2h}$ \\
$u'({\hat{x}}) = 0$\\\\
\end{tabular}
}\hspace{-0.4cm}
\subfloat{%
  \includegraphics[width=0.6\textwidth]{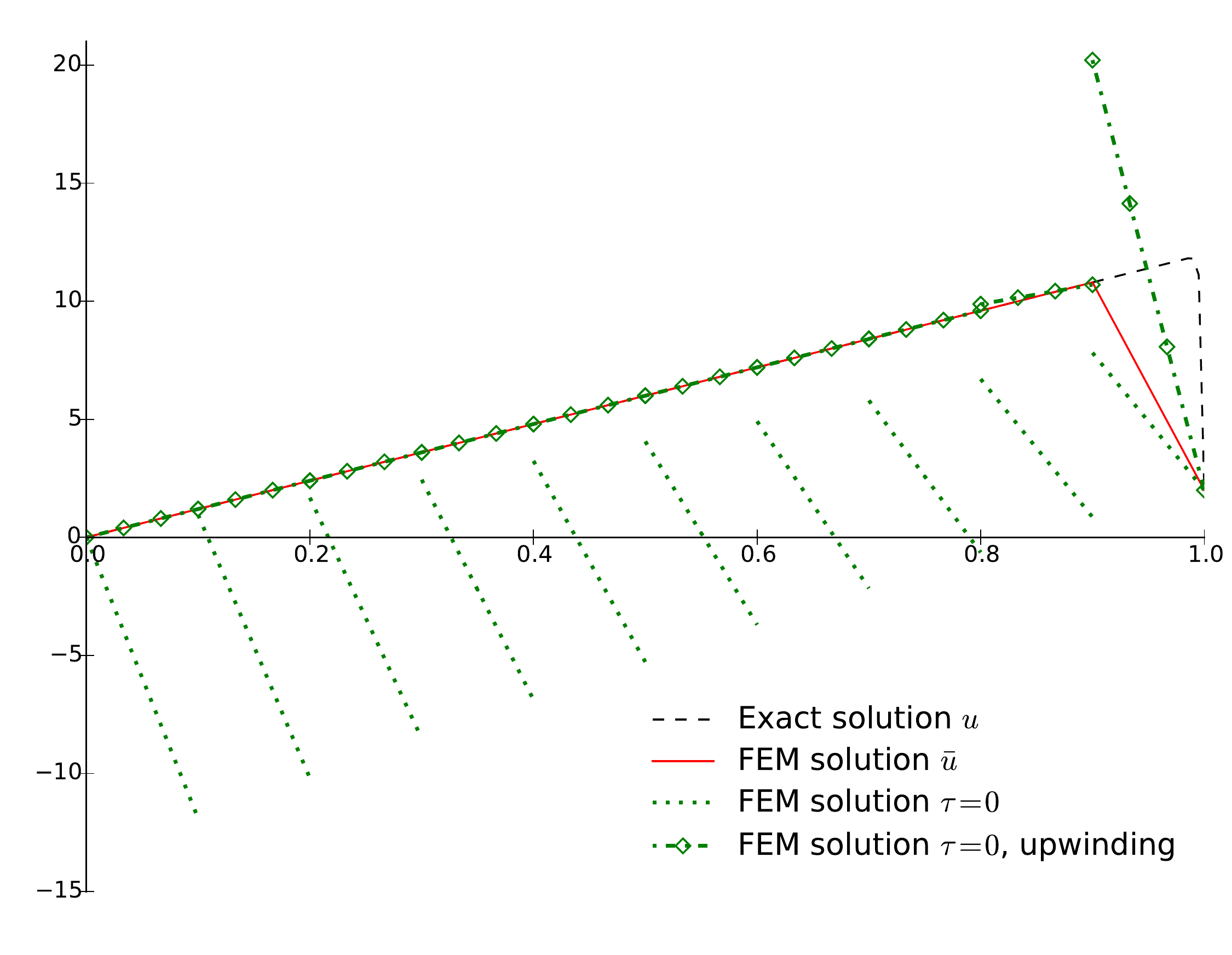}%
}
  \caption{Numerical experiment 3: comparison of the effect of upwinding with the effect of the volumetric fine-scale term $\tau$.}%
  \label{fig:AdDifexp2}
\end{figure}

To illustrate the effect of \eqref{exacttoupwind} on the coarse-scale solution, we consider a third numerical experiment defined in \Cref{fig:AdDifexp2}.
We discretize \cref{coarsescale2AD} with ten linear DG elements, where all fine-scale terms emanating from the diffusion operator are treated as discussed for the $H^1$ projection in \cref{ssec:5.2}. We observe in \Cref{fig:AdDifexp2} that when the fine-scale volumetric term is taken into account, a nodally exact coarse-scale solution $\bar{u}$ is retrieved. When the fine-scale volumetric term is omitted, i.e., $\tau=0$, the coarse-scale solution loses accuracy with respect to the exact solution. 
When an upwind numerical flux is implemented by substitution of \eqref{exacttoupwind},
then the coarse-scale solution is almost indistinguishable from the exact solution in most of the domain, even though $\tau=0$.

The results shown in \Cref{fig:AdDifexp2} indicate that upwind numerical fluxes are able to eliminate the need for a volumetric fine-scale model. 
 This hypothesis may be confirmed by inspecting the residual. In the first nine elements the slope of the coarse-scale solution is such that the residual is nearly zero. The element on the right, however, has a negative slope, whereas $f$ is positive. Therefore, the residual \mbox{$\mathcal{R}_{\bar{u}} = f-a\cdot \nabla \bar{u}$} is large. This means that only in the rightmost element, the fine-scale volumetric term has a significant effect. The upwind method sweeps through the mesh, from left to right. Every subsequent element uses the coarse-scale boundary values of the previous element as boundary conditions. The near-zero residuals on the left mean that the cumulative error due to omission of the volumetric fine-scale term remains small. Upwind numerical fluxes can thus be interpreted as a \textit{remedy} for the lack of the volumetric fine-scale term.



Moreover, upwind numerical fluxes can be interpreted as an implicit fine-scale interface model that can be obtained by manipulating \cref{exacttoupwind} such that the fine-scale term acts as the difference between the upwind term and the average term. We can identify two different cases at an element interface: $a$ points into the element denoted by $K^+$ (where \mbox{$a\cdot n^+ < 0$}) or $a$ points out of that element and into $K^-$ (and thus~\mbox{$a\cdot n^+ > 0$}).\\

{\centering
\begin{tabular}{c|c}
$ \begin{aligned} a\cdot n^+ < 0 \end{aligned}$ & $ \begin{aligned} a\cdot n^+ > 0 \end{aligned}$\\[0.1cm]\hline\hline\\[-0.45cm]
\small$ \begin{aligned} \Bigg. \big\langle\jump{\bar{w}}\cdot a,\bar{u}(\limeps\,\, x-\epsilon a)-\avg{\bar{u}} \big\rangle_{\Gamma}  \end{aligned}= \quad$ &
\small$ \begin{aligned} \big\langle\jump{\bar{w}}\cdot a,\bar{u}(\limeps\,\, x-\epsilon a)-\avg{\bar{u}} \big\rangle_{\Gamma}  \end{aligned}= $ \\[-0.2cm]
\small$  \begin{aligned}  \big\langle\jump{\bar{w}}\cdot a,\bar{u}^- - \frac{1}{2}(\bar{u}^+ + \bar{u}^-)\big\rangle_{\Gamma}\end{aligned}= \quad$ & 
\small$ \quad \begin{aligned} \big\langle\jump{\bar{w}}\cdot a,\bar{u}^+ -\frac{1}{2}(\bar{u}^+ + \bar{u}^- ) \big\rangle_{\Gamma}\end{aligned}=$  \\[0.1cm]
\small$ \begin{aligned} \frac{1}{2}\big\langle\jump{\bar{w}}\cdot a,-\bar{u}^+ + \bar{u}^- \big\rangle_{\Gamma}\end{aligned}= \quad$ & 
\small$ \quad \begin{aligned} \frac{1}{2} \big\langle\jump{\bar{w}}\cdot a,\bar{u}^+ - \bar{u}^- \big\rangle_{\Gamma}\end{aligned}=$  \\[0.1cm]
\small$ \begin{aligned} \frac{1}{2} \big\langle\jump{\bar{w}}\cdot a,n^- \cdot n^+\bar{u}^+ + n^- \cdot n^- \bar{u}^- \big\rangle_{\Gamma}\end{aligned}= \quad$ & 
\small$ \quad \begin{aligned} \frac{1}{2} \big\langle\jump{\bar{w}}\cdot a,n^+ \cdot n^+\bar{u}^+ + n^+ \cdot n^- \bar{u}^- \big\rangle_{\Gamma}\end{aligned}=$  \\[0.0cm]
\small$ \begin{aligned} \frac{1}{2}\big\langle\jump{\bar{w}},\jump{\bar{u}} \, a\cdot n^- \big\rangle_{\Gamma}\end{aligned}= \quad$ & 
\small$ \begin{aligned} \quad \frac{1}{2}\big\langle\jump{\bar{w}},\jump{\bar{u}} \, a\cdot n^+ \big\rangle_{\Gamma}\end{aligned}=$ \\[-0.2cm] 
\end{tabular}
\begin{align}
\frac{1}{2}\big\langle\jump{\bar{w}},\jump{\bar{u}} \, | a\cdot n | \big\rangle_{\Gamma} \label{upwFWM}
\end{align}}

The result of the previous derivation demonstrates that, in either case, the fine-scale model can be retrieved from the following identity:
\begin{align}
\quad\begin{alignedat}{3}
\big\langle\jump{\bar{w}}\cdot a,\avg{u'}\big\rangle_{\Gamma} = \frac{1}{2}\big\langle\jump{\bar{w}},\jump{\bar{u}} \, | a\cdot n | \big\rangle_{\Gamma}\hspace{2.26cm} \\ 
\big\langle\jump{\bar{w}}, (u'^+ + u'^- )\,a  \big\rangle_{\Gamma} = -\big\langle\jump{\bar{w}}, (u'^+ n^+ + u'^-n^- ) | a\cdot n | \big\rangle_{\Gamma} &\\
\big\langle\jump{\bar{w}}, a u'^+ + a u'^-  +| a\cdot n |\, u'^+ n^+ +|\, a\cdot n |  u'^-n^-  \big\rangle_{\Gamma}= 0 \,\,\,&\\
\big\langle\jump{\bar{w}}, a\cdot n^+ \, u'^+ -  a\cdot n^-  u'^-  +| a\cdot n |\, u'^+  -| a\cdot n | \, u'^- \big\rangle_{\Gamma}= 0 \,\,\,&\\
\big\langle\jump{\bar{w}}, (a\cdot n^+ +| a\cdot n |)\, u'^+ - (a\cdot n^- +| a\cdot n | )\, u'^- \big\rangle_{\Gamma}= 0\,\,\, &
\end{alignedat}
\end{align}

When $a$ points out of $K^-$ and into $K^+$, then \mbox{$a\cdot n^+ < 0$} and \mbox{$a\cdot n^+ +| a\cdot n |=0$} which means that the fine-scale component $u'^+$ is removed from the identity. Instead, if $a$ points into $K^-$ and out of $K^+$, then $a\cdot n^- +| a\cdot n |=0$ and the fine-scale component $u'^-$ vanishes. These observations motivate the following the fine-scale model:
\begin{align}
\quad \begin{cases}
\,\, \big\langle\jump{\bar{w}},\, u'^+ \big\rangle_{\Gamma}= 0 \quad\Rightarrow\quad \big\langle\jump{\bar{w}},\, \bar{u}^+ \big\rangle_{\Gamma}=\big\langle\jump{\bar{w}},\, u \big\rangle_{\Gamma} \quad\quad &\text{when } a\cdot n^+ > 0 \\
\,\, \big\langle\jump{\bar{w}},\, u'^- \big\rangle_{\Gamma}= 0 \quad\Rightarrow\quad \big\langle\jump{\bar{w}},\, \bar{u}^- \big\rangle_{\Gamma}=\big\langle\jump{\bar{w}},\, u \big\rangle_{\Gamma} \quad\quad &\text{when } a\cdot n^- > 0 
\end{cases}\label{UpwindModel}
\end{align}
This model only enforces a condition on $u'$ on parts of the element boundary only, namely at all locations where the velocity vector points out of the element. Wherever it points inwards, no condition is imposed on the fine-scale solution. This multiscale interpretation of upwind flux evaluation is based on the assumption that a suitable fine-scale volumetric model is used.


\subsection{Multiscale interpretation of a combined interior penalty method and upwinding formulation}
\label{ssec:5.4}

Classical variational formulations for the advection-diffusion problem \eqref{AdDifProb} rely on both upwinding and proper treatment of the second-order elliptic operator. In this subsection, we investigate the fine-scale interface models that result from using upwind fluxes in combination with the interior penalty method. Our analysis assumes that the volumetric fine-scale term is treated appropriately, according to \cref{taugammagamma2}.



To identify fine-scale models, we compare the variational coarse-scale formulation \eqref{coarsescale2AD} with the classical interior penalty formulation with standard upwind fluxes. After equating the two formulations and eliminating equivalent terms, we can identify the following relation:
\begin{align}
\quad\begin{split}
&\big\langle \jump{\nabla \bar{w}}, \nu \avg{u' } \big\rangle_{\Gamma}+ \big\langle \jump{\bar{w}},a\avg{\bar{u}} \big\rangle_{\Gamma} + \big\langle \jump{\bar{w}},a\avg{u'} \big\rangle_{\Gamma} - \big\langle \jump{\bar{w}},\nu \avg{\nabla u'} \big\rangle_{\Gamma} =\\
&\hspace{2.8cm}\big\langle \jump{\bar{w}},a\,\bar{u}(\limeps\,\, x-\epsilon a)\big\rangle_{\Gamma} + \big\langle \jump{\bar{w}},\nu \eta h^{-1} \jump{\bar{u}} \big\rangle_{\Gamma} \quad\forall\, \bar{w} \in \widebar{\mathcal{V}}(0)
\end{split}
\end{align}
By moving all terms to the left-hand side, and collecting terms, we obtain: 
\begin{align}
\quad\big\langle \jump{\nabla \bar{w}}, \nu \avg{u' } \big\rangle_{\Gamma}+\big\langle \jump{\bar{w}}, a\avg{u'} - \frac{1}{2} \jump{\bar{u}} \, | a\cdot n |  - \nu \avg{\nabla u'} -\nu \eta h^{-1} \jump{\bar{u}} \big\rangle_{\Gamma} = 0
\end{align}
where the upwind terms simplify as in \cref{upwFWM}.

By restricting ourselves again to the one-dimensional case we can make an argument similar to that of \cref{ssec:4.3}. In practice this means that variational terms that involve a different test function operator need to be zero individually. This leaves the following identities at element interfaces, representing the fine-scale interface model:
\begin{align}
\begin{cases}
\,\, \avg{u'}\Big|_{\hat{x}} = 0 \quad & \hat{x} \in \Gamma\\[0.2cm]
\,\, \nu \avg{\nabla u'}\Big|_{\hat{x}} = - \big( \frac{1}{2}  | a\cdot n | + \nu \eta h^{-1} \big)\jump{\bar{u}}\Big|_{\hat{x}}  \qquad & \hat{x} \in \Gamma
\end{cases} \label{FMS_IPup}
\end{align}
We define a distance $d$ as:
\begin{align}
d = \frac{h}{ |a\cdot n|\nu^{-1} h + 2 \eta } \label{FMS_IPupd}
\end{align}
and manipulate the second line in \cref{FMS_IPup} as follows:
\begin{align}
\begin{split}
\avg{\nabla u'}\Big|_{\hat{x}} &= \frac{1}{2}d^{-1}\jump{u'}\Big|_{\hat{x}}\\
-\frac{1}{2}d^{-1} u'^+ n^+ + \frac{1}{2}  \nabla u'^+    &=  \frac{1}{2}d^{-1} u'^- n^- - \frac{1}{2} \nabla u'^-\\
u'^{+} - d\,\nabla u'^+ n^{+} &=u'^{-} - d\, \nabla u'^- n^{-} \\
u'^+(\hat{x}-d\, n^{+} ) &\approx u'^-(\hat{x}-d\, n^{-} ) \label{FMS_IPupImp}
\end{split}
\end{align}

We observe that the fine-scale interface model \eqref{FMS_IPup} can be interpreted in exactly the same way as the fine-scale model \cref{IPmodel} of the IP method for the Poisson problem. The influence of the advection operator manifests itself through the distance $d$ that includes the ratio between the diffuse parameter $\nu$ and the velocity $a$. It is obvious that increasing velocity at constant $\nu$ has an additional ``clamping'' effect. 
The definition \eqref{FMS_IPupd} correctly reduces to the prior expression \cref{dIP}, when the velocity $a$ is zero. 

\begin{figure}[ht]
\hspace{-0.2cm}\subfloat{%
\begin{tabular}{c}\vspace{-7cm}\\
\textbf{Model problem} \\
$f=6$\\
$a=-0.5$\\
$\nu=0.15$ \\
$x_0 = 0$, $x_1 = 0.9$\\
$u_0 = 0$, $u_1 = 2$ \\
$d = 0.1$\\\\
\textbf{Exact solution} \\
$u(x) =$  \\
$\dfrac{u_1 \!-\! \frac{f}{a}x_1}{e^{\frac{a}{\nu}x_1}\!-\! 1}\left(e^{\frac{a}{\nu}x} \!-\! 1\right)+\frac{f}{a}x$  
\end{tabular}
}\hspace{-0.2cm}
\subfloat{%
  \includegraphics[width=0.7\textwidth]{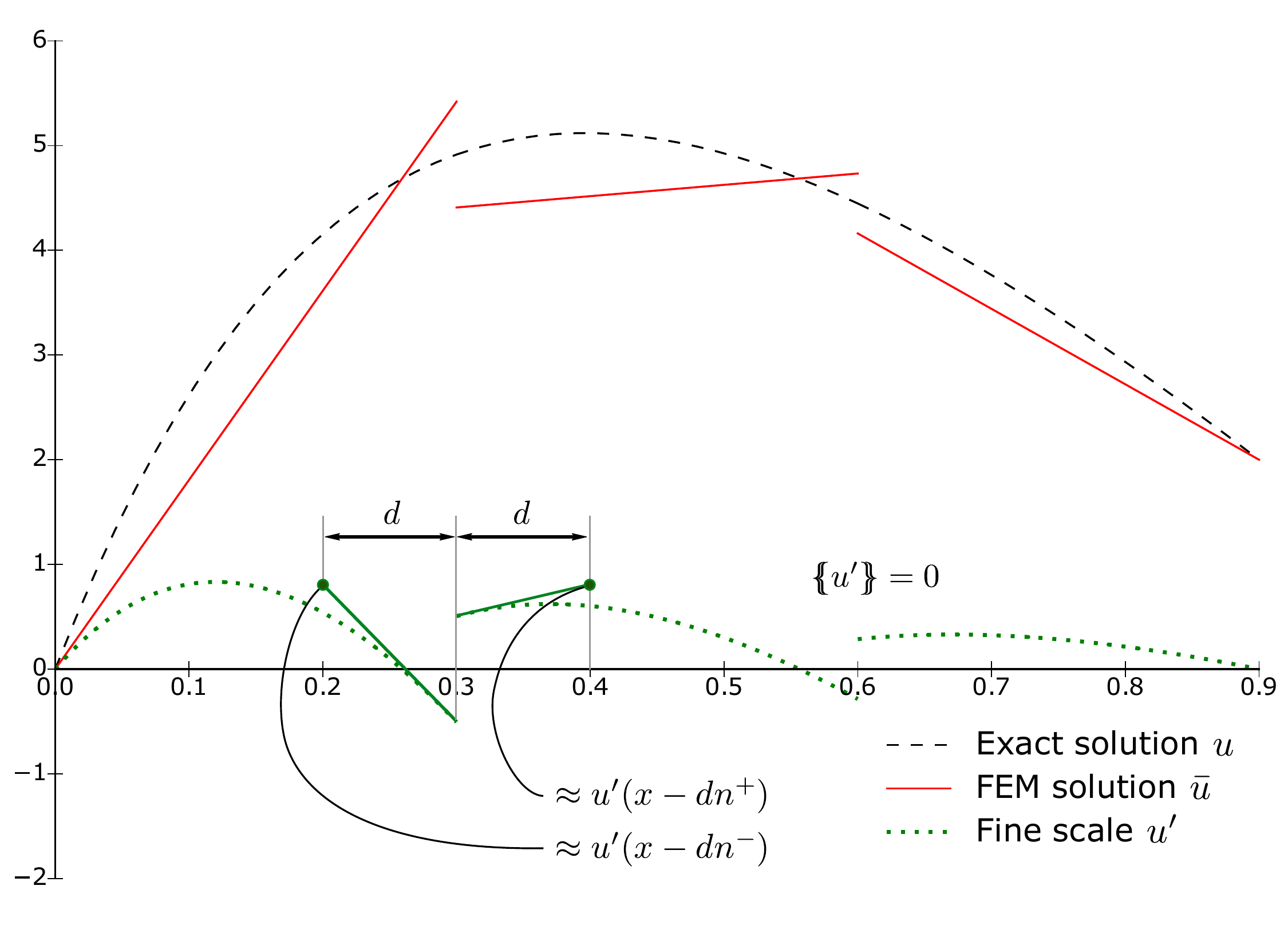}%
}
  \caption{IP method with upwinding: Coarse- and fine-scale solutions of a three-element linear DG discretization and graphical interpretation of the effect of the associated fine-scale interface model.}%
  \label{fig:AdDif_IP}
\end{figure}

The model problem defined in \Cref{fig:AdDif_IP} illustrates the effect of the fine-scale interface model \eqref{FMS_IPup} that corresponds to the IP method with upwinding. Using a three element linear DG discretization and a constant forcing term, we obtain the solution of the variational coarse-scale formulation \eqref{coarsescale2AD} with \eqref{FMS_IPup}. We reiterate that we treat the fine-scale volumetric term implicitly with relation \cref{taugammagamma2}. The coarse- and fine-scale solutions in \Cref{fig:AdDif_IP} confirm that the fine-scale solution satisfies \cref{FMS_IPup,FMS_IPupImp}.


\section{Summary and conclusions}
\label{sec:conclusions}

In this article, we have developed a general strategy for obtaining VMS formulations suitable for discontinuous Galerkin discretizations. We transferred the original VMS idea into a DG setting, that is, the decomposition of the true solution into a discontinuous coarse-scale solution and an accompanying discontinuous fine-scale solution. We defined the associated weak formulations on a per-element basis, and used multiscale-type transmission conditions to couple the elements.  

We obtained coarse-scale weak formulations that include fine-scale volumetric terms and fine-scale interface terms. These fine-scale contributions are defined by the projector that was used to decompose the true solution into coarse-scale and fine-scale components. We extended the existing residual-based model of the fine-scale volumetric term to include nonhomogeneous fine-scale element boundary values. We also replaced the fine-scale interface terms by explicit interface models (known values) or implicit interface models (relations with coarse-scale components). 

In our numerical experiments we showed that these models form a set of conditions to which the fine-scale solution must adhere. In particular, for the one-dimensional Poisson problem the fine-scale interface models become pointwise conditions. For multi-dimensional problems the models combine into a single more complex condition that is satisfied weakly. In the future, we will investigate the multi-dimensional variant of the fine-scale condition further.



An important side-effect of our multiscale formulation is that it naturally enables new perspectives on discontinuous Galerkin finite element formulations and their numerical properties. In particular, we rederived the symmetric interior penalty method for a Poisson problem from a multiscale viewpoint by choosing a particular multiscale projector. 
Additionally, we showed that the use of upwind fluxes in advection-diffusion problems helps to compensate for the absence of a suitable fine-scale volumetric model. This observation could be a new point of view for explaining the effectiveness of upwind formulations. 





The multiscale DG framework developed in this work can be the foundation for further research in a number of directions.
For example, new DG formulations could be developed by constructing fine-scale interface models that retrieve the coarse-scale solution corresponding to some favorable projection.
Furthermore, the effectiveness of the residual-based volumetric fine-scale model as a multiscale model may be evaluated, in particular with respect to the nonhomogeneous fine-scale boundary values introduced in this work.
Another promising idea is the extension of the theory to a mixed method setting, as many advanced DG formulations for the Poisson problem are derived from this point of view. This could enable the derivation of additional formulations and provide more flexibility in manipulating the fine-scale terms.

\appendix
\section{Green's function for the 1-D advection-diffusion problem}
\label{App:A}
In this appendix, we provide the full details on the derivation of the Green's function associated with the advection-diffusion problem of \Cref{sec:5} and obtain analytical expressions for the corresponding parameters $\tau$, $\gamma_0$ and $\gamma_1$. 

The advection-diffusion equation involves the following differential operator, adjoint operator, and accompanying interface terms.
\begin{align}
\begin{split}
&\hspace{2cm}\mathcal{L} = a\cdot\nabla -\nu\Delta; \qquad \mathcal{L}^* = -a\cdot\nabla -\nu\Delta;   \\
&k( w',u';\partial K) = \big\langle w', a\cdot n u'\big\rangle_{\partial K} -\big\langle w', \nu n\cdot u'_{y} \big\rangle_{\partial K} + \big\langle n\cdot \dd{y} w' , \nu u' \big\rangle_{\partial K}\label{app:addif}
\end{split}
\end{align}
Also recall the definition of the Green's function:
\begin{align}
\begin{cases}
\, g(x,y) \in \mathcal{V}(0)\\
\, \mathcal{L}^*\! g(x,y) = \delta_x \quad&\text{for } y\in K\\
\, g(x,y) = 0  \quad \qquad &\text{for } y\in \partial K
\end{cases}\label{gDef2}
\end{align}

When we substitute the Green's function in place of the test function in the term $k( w',u';\partial K)$ defined in \cref{app:addif}, it reduces to:
\begin{align}
k( g(x,y) ,u';\partial K) = \big\langle n\cdot \dd{y} g(x,y) , \nu u' \big\rangle_{\partial K}
\end{align}
This is in line with the definition of $\gamma$ in \cref{gammadef}.

We can then split the Green's function into two components:
\begin{align}
g(x,y) = \begin{cases}\,\,\, g_1 (x,y)  \qquad \text{when } y<x \\
\,\,\, g_2 (x,y)  \qquad \text{when } y\geq x \end{cases}
\end{align}
By definition, $\mathcal{L}^*\! g(x,y) = \delta_x$, and therefore, the following holds:
\begin{align}
\mathcal{L}^*\! g_1 (x,y) = -a \dd{y} g_1 (x,y) - \nu \ddd{y}  g_1 (x,y)  =0, \qquad \text{when } y<x \\
\mathcal{L}^*\! g_2 (x,y) = -a \dd{y} g_2 (x,y) - \nu \ddd{y}  g_2 (x,y)  =0,  \qquad \text{when } y> x
\end{align}
Both sections of the Green's functions are thus of the form:
\begin{align}
g_1 (x,y) = C_1 e^{-\frac{a}{\nu} y}+D_1 \\
g_2 (x,y) = C_2 e^{-\frac{a}{\nu} y}+D_2
\end{align}\\[-0.8cm]

\noindent \textbf{Condition 1 \& 2:} The boundary conditions in \cref{gDef2} set the Green's function to zero on the element boundary. They can be used to determine two unknown coefficients:
\begin{alignat}{5}
&g_1 (x,x_j) &&= C_1 e^{-\frac{a}{\nu} x_j}+D_1 &&= 0 \quad \Rightarrow \quad D_1 = - C_1 e^{-\frac{a}{\nu} x_j}\\
&g_2 (x,x_{j+1}) &&= C_2 e^{-\frac{a}{\nu} x_{j+1}}+D_2 &&= 0 \quad \Rightarrow \quad D_2 = -C_2 e^{-\frac{a}{\nu} x_{j+1}}
\end{alignat}\\[-0.8cm]

\noindent \textbf{Condition 3:} The Green's function must be continuous. Otherwise, the second derivatives in the differential operator would lead to derivatives of Dirac $\delta$ distributions, rather than just Dirac $\delta$ distributions. Notice that this condition also ensures that $g(x,y) \in \mathcal{V}(0)$, where $\mathcal{V}(0)$ is defined in \cref{V}. Continuity of $g(x,y)$ requires:
\begin{align}
&g_1 (x,x) = C_1\left(  e^{-\frac{a}{\nu} x} - e^{-\frac{a}{\nu} x_j}\right) = g_2 (x,x) = C_2\left( e^{-\frac{a}{\nu} x}- e^{-\frac{a}{\nu} x_{j+1}} \right)\\
&\quad \Rightarrow C_2 = C_1  \frac{e^{-\frac{a}{\nu} x} - e^{-\frac{a}{\nu} x_j}}{e^{-\frac{a}{\nu} x}- e^{-\frac{a}{\nu} x_{j+1}}} \label{req1}
\end{align}

\noindent \textbf{Condition 4:} The final condition is obtained from the following Dirac delta distribution property:
\begin{align}
\int\limits_{x\uparrow}^{x\downarrow} \mathcal{L}^*\! g(x,y) \d{y} = \int\limits_{x\uparrow}^{x\downarrow} \delta_x = 1
\end{align}
\begin{align}
-a g (x,y) \Big|_{x\uparrow}^{x\downarrow} - \nu \dd{y}  g (x,y)  \Big|_{x\uparrow}^{x\downarrow} = \nu \dd{y}  g_1 (x,y)\Big|_x - \nu \dd{y}  g_2 (x,y)\Big|_x = 1 \label{eqGdirac}
\end{align}
The first term on the left in \cref{eqGdirac} cancels due to the continuity condition. This results in:
\begin{align}
-a C_1 e^{-\frac{a}{\nu} x} + a C_2 e^{-\frac{a}{\nu} x}  = 1 \label{req2}
\end{align}
Combining relations \eqref{req1} and \eqref{req2} yields the following expression for the Green's function:
\begin{align}
g(x,y) = \begin{cases}\,\,\, g_1(x,y) = \dfrac{b}{k-1} \left( e^{-\frac{a}{\nu} y}- e^{-\frac{a}{\nu} x_{j}} \right)  \qquad &\text{when } y<x \\
\,\,\, g_2(x,y) = \dfrac{k\,b}{k-1} \left(  e^{-\frac{a}{\nu} y} - e^{-\frac{a}{\nu} x_{j+1}}\right)  \qquad &\text{when } y\geq x \end{cases}
\end{align}
where:
\begin{align} 
&k = \frac{e^{-\frac{a}{\nu} x} - e^{-\frac{a}{\nu} x_j}}{e^{-\frac{a}{\nu} x}- e^{-\frac{a}{\nu} x_{j+1}}} \\[0.2cm]
&b = \frac{1}{a} e^{\frac{a}{\nu}x}
\end{align}

Finally, this Green's function is used to obtain the element constants $\tau$, $\gamma_0$ and $\gamma_1$.  A single expression for these integral quantities can be obtained by translating the element such that $x_j = 0$ and $x_{j+1}=h$. With the help of definitions \eqref{taudef} and \eqref{gammadef}, we obtain the following explicit expressions:
\begin{align}
\small\begin{alignedat}{5}
&\tau\, &&= \frac{1}{h} \int\limits_{0}^h \!\!\int\limits_{0}^h g(x,y) \dx\d{y}  &&= \frac{1}{h} \int\limits_{0}^h \!\left[\int\limits_{0}^x g_1(x,y) \d{y} \!+\! \!\int\limits_{x}^h g_2(x,y) \d{y} \right]\!\!\dx = \frac{h}{2\,a}\!-\!\frac{\nu}{a^2}\!+\!\frac{h}{a \left( e^{\frac{a}{\nu}h}\!-\!1 \right)} \label{apptauADD}\\
&\gamma_0 &&= \frac{1}{h} \int\limits_{0}^h \dd{y} g(x,0)\dx  &&= \frac{1}{h} \int\limits_{0}^h \dd{y} g_{1}(x,0) \dx = - \frac{1}{h} \int\limits_{0}^h  \frac{a}{\nu} \frac{b}{k-1}\dx   =  \frac{\nu - a\, h - \nu e^{-\frac{a}{\nu} h}}{a\, h\, \nu \left(e^{-\frac{a}{\nu} h}-1\right)}  \\
&\gamma_1 &&= \frac{1}{h} \int\limits_{0}^h \dd{y} g(x,h)\dx &&= \frac{1}{h} \int\limits_{0}^h \dd{y} g_{2}(x,h) \dx = - \frac{1}{h} \int\limits_{0}^h  \frac{a}{\nu} \frac{k\,b}{k-1}e^{-\frac{a}{\nu} h}\dx =  \frac{\nu + a\, h - \nu e^{\frac{a}{\nu} h}}{a\, h\, \nu \left(e^{\frac{a}{\nu} h}-1\right)} 
\end{alignedat}
\end{align}

\section*{Acknowledgments}
D. Schillinger gratefully acknowledges support from the National Science Foundation through the NSF CAREER Award No. 1651577. The authors thank Prof.\ Henryk Stolarski for his help with the derivation of the Green's function for the advection-diffusion problem.

\newpage

\bibliographystyle{siamplain}
\bibliography{MyBib}

\begin{thebibliography}{10}

\bibitem{Akkerman2008}
{\sc I.~Akkerman, Y.~Bazilevs, V.~M. Calo, T.~J.~R. Hughes, and S.~J.
  Hulshoff}, {\em {The role of continuity in residual--based variational
  multiscale modeling of turbulence}}, Computational Mechanics, 41 (2008),
  pp.~371--378.

\bibitem{Arnold2000}
{\sc D.~Arnold, F.~Brezzi, B.~Cockburn, and D.~Marini}, {\em {Discontinuous
  Galerkin methods for elliptic problems}}, in Discontinuous Galerkin Methods.
  Theory, Computation and Applications, B.~Cockburn, G.~E. Karniadakis, and
  C.~W. Shu, eds., vol.~11 of Lecture Notes in Computational Science and
  Engineering, Springer-Verlag, Berlin, Heidelberg, 2000, pp.~89--101.

\bibitem{Arnold2001}
{\sc D.~N. Arnold, F.~Brezzi, B.~Cockburn, and L.~D. Marini}, {\em {Unified
  Analysis of Discontinuous Galerkin Methods for Elliptic Problems}}, SIAM
  Journal on Numerical Analysis, 39 (2002), pp.~1749--1779.

\bibitem{Babuska1973}
{\sc I.~Babu\v{s}ka and M.~Zl\'amal}, {\em {Nonconforming elements in the
  finite element method with penalty}}, SIAM Journal of Numerical Analysis, 10
  (1973), pp.~863--875.

\bibitem{Bassi1997}
{\sc F.~Bassi and S.~Rebay}, {\em {A High-Order Accurate Discontinuous Finite
  Element Method for the Numerical Solution of the Compressible Navier--Stokes
  Equations}}, Journal of Computational Physics, 131 (1997), pp.~267 -- 279.

\bibitem{Bassi1997b}
{\sc F.~Bassi, S.~Rebay, G.~Mariotti, S.~Pedinotti, and M.~Savini}, {\em {A
  high-order accurate discontinuous finite element method for inviscid and
  viscous turbomachinery flows}}, in Proceedings of 2nd European Conference on
  Turbomachinery, Fluid Dynamics and Thermodynamicst, R.~Decuypere and
  G.~Dibelius, eds., Antwerpen, Belgium, 1997, pp.~99--108.

\bibitem{Baumann1999}
{\sc C.~E. Baumann and J.~T. Oden}, {\em {A discontinuous \textit{hp}-finite
  element method for convection-diffusion problems}}, Computer Methods in
  Applied Mechanics and Engineering, 175 (1999), pp.~311--341.

\bibitem{Bazilevs2006}
{\sc Y.~Bazilevs}, {\em {Isogeometric analysis of turbulence and
  fluid-structure interaction}}, PhD thesis, The University of Texas at Austin,
  2006.

\bibitem{Bazilevs2007}
{\sc Y.~Bazilevs, V.~M. Calo, J.~A. Cottrell, T.~J.~R. Hughes, A.~Reali, and
  G.~Scovazzi}, {\em {Variational multiscale residual-based turbulence modeling
  for large eddy simulation of incompressible flows}}, Computer Methods in
  Applied Mechanics and Engineering, 197 (2007), pp.~173--201.

\bibitem{Bochev2005}
{\sc P.~Bochev, T.~J.~R. Hughes, and G.~Scovazzi}, {\em {A multiscale
  discontinuous Galerkin method}}, in International Conference on Large-Scale
  Scientific Computing, Springer, 2005, pp.~84--93.

\bibitem{Brezzi1997b}
{\sc F.~Brezzi, L.~P. Franca, T.~J.~R. Hughes, and A.~Russo}, {\em $ b=\int g
  $.}, Computer Methods in Applied Mechanics and Engineering,  (1997),
  pp.~329--339.

\bibitem{Brezzi1997}
{\sc F.~Brezzi, G.~Manzini, D.~Marini, P.~Pietra, and A.~Russo}, {\em
  {Discontinuous finite elements for diffusion problems}}, in Atti Convegno in
  onore, F.~Brioschi, ed., Milan, Italy, 1997, pp.~197--217.

\bibitem{Brezzi2000}
{\sc F.~Brezzi, G.~Manzini, D.~Marini, P.~Pietra, and A.~Russo}, {\em
  {Discontinuous Galerkin approximations for elliptic problems}}, Numerical
  Methods for Partial Differential Equations, 16 (2000), pp.~365--378.

\bibitem{Brooks1982}
{\sc A.~N. Brooks and T.~J.~R. Hughes}, {\em {Streamline
  upwind/Petrov--Galerkin formulations for convection dominated flows with
  particular emphasis on the incompressible Navier--Stokes equations}},
  Computer Methods in Applied Mechanics and Engineering, 32 (1982),
  pp.~199--259.

\bibitem{Buffa2006}
{\sc A.~Buffa, T.~J.~R. Hughes, and G.~Sangalli}, {\em {Analysis of a
  multiscale discontinuous Galerkin method for convection-diffusion problems}},
  SIAM Journal on Numerical Analysis, 44 (2006), pp.~1420--1440.

\bibitem{Calo2004}
{\sc V.~M. Calo}, {\em {Residual--based Multiscale Turbulence modeling: Finite
  Volume simulations of bypass transition}}, PhD thesis, Stanford University,
  2004.

\bibitem{Cockburn1999}
{\sc B.~Cockburn, E.~K. George, and C.-W. Shu}, {\em {Discontinuous Galerkin
  methods: Theory, Computation and Applications}}, Springer, Berlin,
  Heidelberg, 11th~ed., 2000.

\bibitem{Cockburn2009}
{\sc B.~Cockburn, J.~Gopalakrishnan, and R.~Lazarov}, {\em {Unified
  hybridization of discontinuous Galerkin, mixed, and continuous Galerkin
  methods for second order elliptic problems}}, SIAM Journal on Numerical
  Analysis, 47 (2009), pp.~1319--1365.

\bibitem{Cockburn:05.1}
{\sc B.~Cockburn, G.~Kanschat, and D.~Sch{\"o}tzau}, {\em The local
  discontinuous {G}alerkin method for linearized incompressible fluid flow: a
  review}, Computers \& Fluids, 34 (2005), pp.~491--506.

\bibitem{Cockburn:00.1}
{\sc B.~Cockburn, G.~Karniadakis, and C.-W. Shu~(eds.)}, {\em Discontinuous
  {G}alerkin methods: theory, computation and applications}, Springer, 2000.

\bibitem{Cockburn1998}
{\sc B.~Cockburn and C.-W. Shu}, {\em {The Local Discontinuous Galerkin Method
  for Time--Dependent Convection--Diffusion Systems}}, SIAM Journal on
  Numerical Analysis, 35 (1998), pp.~2440--2463.

\bibitem{Coley2017}
{\sc C.~Coley and J.~A. Evans}, {\em {Variational Multiscale Modeling with
  Discontinuous Subscales: Analysis and Application to Scalar Transport}},
  arXiv preprint arXiv:1705.00082,  (2017).

\bibitem{Collis2002}
{\sc S.~S. Collis}, {\em {The DG/VMS Method for Unified Turbulence
  Simulation}}, in 32nd AIAA Fluid Dynamics Conference and Exhibit, Reston,
  Virigina, 2002.

\bibitem{Collis2015}
{\sc S.~S. Collis and S.~Ramakrishnan}, {\em {The Local Variational Multiscale
  Method}}, Third MIT Conference on Computational Fluid and Solid Dynamics,
  (2005).

\bibitem{Donea2003}
{\sc J.~Don\'ea and A.~Huerta}, {\em {Finite element methods for flow
  problems}}, John Wiley \& Sons, Ltd, Hoboken, New Jersey, 2003.

\bibitem{Douglas1976}
{\sc J.~Douglas and T.~Dupont}, {\em {Interior penalty procedures for elliptic
  and parabolic Galerkin methods}}, Springer, Berlin, Heidelberg, 1976,
  pp.~207--216.

\bibitem{Hartmann:08.1}
{\sc R.~Hartmann}, {\em Numerical analysis of higher order {D}iscontinuous
  {G}alerkin finite element methods}, in {CFD - ADIGMA} course on very high
  order discretization methods, H.~Deconinck, ed., vol.~VKI LS 2008-08 of Von
  Karman Institute for Fluid Dynamics, Belgium, 2008, pp.~1--107.

\bibitem{Hsu2012a}
{\sc M.-C. Hsu and Y.~Bazilevs}, {\em {Fluid--structure interaction modeling of
  wind turbines: Simulating the full machine}}, Computational Mechanics, 50
  (2012), pp.~821--833.

\bibitem{Hsu2014}
{\sc M.-C. Hsu, D.~Kamensky, Y.~Bazilevs, M.~S. Sacks, and T.~J.~R. Hughes},
  {\em {Fluid--structure interaction analysis of bioprosthetic heart valves:
  significance of arterial wall deformation}}, {Computational mechanics}, 54
  (2014), pp.~1055--1071.

\bibitem{Huerta2013}
{\sc A.~Huerta, A.~Angeloski, X.~Roca, and J.~Peraire}, {\em {Efficiency of
  high-order elements for continuous and discontinuous Galerkin methods}},
  International Journal for Numerical Methods in Engineering, 96 (2013),
  pp.~529--560.

\bibitem{Hughes2003}
{\sc T.~J. Hughes and A.~A. Oberai}, {\em {Calculation of shear stresses in the
  Fourier-Galerkin formulation of turbulent channel flows: projection, the
  Dirichlet filter and conservation}}, Journal of Computational Physics, 188
  (2003), pp.~281--295.

\bibitem{Hughes1995}
{\sc T.~J.~R. {Hughes}}, {\em {Multiscale phenomena: Green's functions, the
  Dirichlet--to--Neumann formulation, subgrid scale models, bubbles and the
  origins of stabilized methods}}, Computer Methods in Applied Mechanics and
  Engineering, 127 (1995), pp.~387--401.

\bibitem{Hughes1998}
{\sc T.~J.~R. Hughes, G.~R. Feij\'oo, L.~{Mazzei}, and J.-B. Quincy}, {\em {The
  variational multiscale method -- a paradigm for computational mechanics}},
  Computer Methods in Applied Mechanics and Engineering, 166 (1998), pp.~3--24.

\bibitem{HUGHES1989173}
{\sc T.~J.~R. Hughes, L.~P. Franca, and G.~M. Hulbert}, {\em {A new finite
  element formulation for computational fluid dynamics: VIII. The
  Galerkin/least--squares method for advective--diffusive equations}}, Computer
  Methods in Applied Mechanics and Engineering, 73 (1989), pp.~173--189.

\bibitem{Hughes1999}
{\sc T.~J.~R. Hughes, L.~Mazzei, and K.~E. Jansen}, {\em {Large Eddy Simulation
  and the variational multiscale method}}, Computing and Visualization in
  Science, 3 (2000), pp.~47--59.

\bibitem{Hughes2001a}
{\sc T.~J.~R. Hughes, L.~Mazzei, A.~A. Oberai, and A.~A. Wray}, {\em {The
  multiscale formulation of large eddy simulation: Decay of homogeneous
  isotropic turbulence}}, Physics of Fluids, 13 (2001), pp.~505--512.

\bibitem{Hughes2001}
{\sc T.~J.~R. Hughes, A.~A. Oberai, and L.~Mazzei}, {\em {Large eddy simulation
  of turbulent channel flows by the variational multiscale method}}, Physics of
  Fluids, 13 (2001), pp.~1784--1799.

\bibitem{Hughes2006}
{\sc T.~J.~R. Hughes, G.~Scovazzi, P.~B. Bochev, and A.~Buffa}, {\em {A
  multiscale discontinuous Galerkin method with the computational structure of
  a continuous Galerkin method}}, Computer Methods in Applied Mechanics and
  Engineering, 195 (2006), pp.~2761--2787.

\bibitem{Hughes2004b}
{\sc T.~J.~R. Hughes, G.~Scovazzi, and L.~P. Franca}, {\em Multiscale and
  stabilized methods}, in Encyclopedia of computational mechanics, E.~Stein,
  R.~De~Borst, and T.~J.~R. Hughes, eds., John Wiley \& Sons, Ltd, 2004, ch.~4.

\bibitem{Hughes1996}
{\sc T.~J.~R. Hughes and J.~R. Stewart}, {\em {A space--time formulation for
  multiscale phenomena}}, Journal of Computational and Applied Mathematics, 74
  (1996), pp.~217--229.

\bibitem{Kamensky2015}
{\sc D.~Kamensky, M.-C. Hsu, D.~Schillinger, J.~A. Evans, A.~Aggarwal,
  Y.~Bazilevs, M.~S. Sacks, and T.~J.~R. Hughes}, {\em {An immersogeometric
  variational framework for fluid--structure interaction: Application to
  bioprosthetic heart valves}}, {Computer methods in applied mechanics and
  engineering}, 284 (2015), pp.~1005--1053.

\bibitem{Kirby2012}
{\sc R.~M. Kirby, S.~J. Sherwin, and B.~Cockburn}, {\em {To CG or to HDG: A
  comparative study}}, Journal of Scientific Computing, 51 (2012),
  pp.~183--212.

\bibitem{Lehrenfeld:15.1}
{\sc C.~Lehrenfeld and J.~Sch{\"o}berl}, {\em High order exactly
  divergence-free hybrid discontinuous {G}alerkin methods for unsteady
  incompressible flows}, Computer Methods in Applied Mechanics and Engineering,
  307 (2016), pp.~339--361.

\bibitem{Mavriplis:09.1}
{\sc D.~Mavriplis, C.~Nastase, K.~Shahbazi, L.~Wang, and N.~Burgess}, {\em
  Progress in high-order discontinuous {G}alerkin methods for aerospace
  applications}, AIAA paper, 601 (2009).

\bibitem{Munts}
{\sc E.~A. Munts, S.~J. Hulshoff, and R.~{De Borst}}, {\em {A Space-Time
  Variational Multiscale Discretization for LES}}, in 34th AIAA Aerospace
  Sciences Meeting and Exhibit, Portland, Oregon, 2004.

\bibitem{Peraire2007}
{\sc J.~Peraire and P.-O. Persson}, {\em {The compact discontinuous Galerkin
  (CDG) method for elliptic problems}}, SIAM Journal on Scientific Computing,
  30 (2007), pp.~1806--1824.

\bibitem{Ramakrishnan2004a}
{\sc S.~Ramakrishnan and S.~S. Collis}, {\em {Turbulence Control Simulation
  Using the Variational Multiscale Method}}, AIAA Journal, 42 (2004),
  pp.~745--753.

\bibitem{Reed1973}
{\sc W.~H. Reed and T.~R. Hill}, {\em {Triangular Mesh Methods for the Neutron
  Transport Equation}}, Proceedings of the American Nuclear Society, 836
  (1973), pp.~1--23.

\bibitem{Riviere1999}
{\sc B.~Rivi\`ere, M.~F. Wheeler, and V.~Girault}, {\em {Improved energy
  estimates for interior penalty, constrained and discontinuous Galerkin
  methods for elliptic problems. Part I}}, Computational Geosciences, 3 (1999),
  pp.~337--360.

\bibitem{Sangalli2004}
{\sc G.~Sangalli}, {\em {A discontinuous residual-free bubble method for
  advection-diffusion problems}}, Journal of Engineering Mathematics, 49
  (2004), pp.~149--162.

\bibitem{Schillinger2016}
{\sc D.~Schillinger, I.~Harari, M.-C. Hsu, D.~Kamensky, S.~K.~F. Stoter, Y.~Yu,
  and Y.~Zhao}, {\em {The non-symmetric Nitsche method for the parameter-free
  imposition of weak boundary and coupling conditions in immersed finite
  elements}}, Computer Methods in Applied Mechanics and Engineering, 309
  (2016), pp.~625--652.

\bibitem{Tezduyar1991}
{\sc T.~E. Tezduyar}, {\em {Stabilized Finite Element Formulations for
  Incompressible Flow Computations}}, Advances in Applied Mechanics, 28 (1991),
  pp.~1--44.

\bibitem{Wang:13.2}
{\sc Z.~J. Wang, K.~Fidkowski, R.~Abgrall, F.~Bassi, D.~Caraeni, A.~Cary,
  H.~Deconinck, R.~Hartmann, K.~Hillewaert, H.~T. Huynh, N.~Kroll, G.~May,
  P.-O. Persson, B.~van Leer, and M.~Visbal}, {\em {High-Order CFD Methods:
  Current Status and Perspective}}, International Journal for Numerical Methods
  in Fluids, 72 (2013), pp.~811--845.

\bibitem{Xu2016}
{\sc F.~Xu, D.~Schillinger, D.~Kamensky, V.~Varduhn, C.~Wang, and M.-C. Hsu},
  {\em {The tetrahedral finite cell method for fluids: Immersogeometric
  analysis of turbulent flow around complex geometries}}, Computers {\&}
  Fluids,  (2015).

\bibitem{Xu2009}
{\sc Y.~Xu and C.~W. Shu}, {\em {Local discontinuous Galerkin methods for
  high-order time-dependent partial differential equations}}, Communications in
  Computational Physics, 7 (2010), pp.~1--46.

\end{thebibliography}
\end{document}